\documentclass[10pt]{article}
\usepackage[margin=0.8in]{geometry}
\usepackage[breaklinks=true,colorlinks=true,hypertexnames=false,linktoc=section,pdfencoding=auto,psdextra=true,linkcolor=blue,anchorcolor=black,citecolor=blue,filecolor=black,menucolor=black,runcolor=black,urlcolor=black]{hyperref}
\usepackage{epsfig,amssymb,soul,url,cite,graphicx,epstopdf,amsfonts,amsmath,bm,setspace,comment,siunitx,mathrsfs,multirow,mathtools,array,lipsum,enumitem,ifthen,leftidx,relsize,braket,float,color,titlesec,titling,microtype,caption,algorithmic,tikz,abstract}
\usepackage[normalem]{ulem}
\usepackage[affil-it]{authblk}
\allowdisplaybreaks
\usepackage{mathpazo}
\usepackage[semibold]{sourcesanspro}
\DeclareRobustCommand{\sbseries}{\fontseries{sb}\selectfont}
\DeclareTextFontCommand{\textsb}{\sbseries}
\titleformat*{\section}{\large\sffamily\sbseries}
\titleformat*{\subsection}{\normalsize\sffamily\sbseries}
\titleformat*{\subsubsection}{\normalsize\sffamily\sbseries}
\titleformat*{\paragraph}{\normalsize\sffamily\sbseries}
\usepackage[framemethod=TikZ]{mdframed}
\mdfdefinestyle{MyFrame}{linecolor=black,outerlinewidth=1pt,roundcorner=1pt,innertopmargin=5pt,innerbottommargin=5pt,innerrightmargin=5pt,innerleftmargin=5pt,backgroundcolor=white}

\let\OLDthebibliography\thebibliography
\renewcommand\thebibliography[1]{\OLDthebibliography{#1}\setlength{\parskip}{0pt}
  \setlength{\itemsep}{3.8pt plus 0.5ex}}

\newcolumntype{C}[1]{>{\centering\arraybackslash}p{#1}}

\usepackage[ruled,vlined]{algorithm2e}
\usepackage[most]{tcolorbox}

\SetCommentSty{mycommfont}
\usepackage{tikz}
\newcommand{\tikzc}[2][black,fill=black]{\tikz[baseline=-0.5ex]\draw[#1,radius=#2] (0,0) circle ;}
\makeatletter
\DeclareRobustCommand{\sqcdot}{\mathbin{\mathpalette\morphic@sqcdot\relax}}
\newcommand{\morphic@sqcdot}[2]{%
  \sbox\z@{$\m@th#1\centerdot$}%
  \ht\z@=.33333\ht\z@
  \vcenter{\box\z@}%
}
\makeatother

\makeatletter
\def\thanks#1{\protected@xdef\@thanks{\@thanks
		\protect\footnotetext{#1}}}
\makeatother

\title{\begin{mdframed}[style=MyFrame]\begin{center}
Multiscale multimesh finite element method \big| $\text{M}^2$-FEM:\\Hierarchical mesh-decoupling for integral structural theories
\end{center}\end{mdframed}}
\author{Wei Ding$^*$\thanks{$^*$To whom correspondence should be addressed. Email: ding242@purdue.edu or fsemperl@purdue.edu}}
\author{Sansit Patnaik}
\author{Fabio Semperlotti$^*$}
\affil{Ray W. Herrick Laboratories, School of Mechanical Engineering\\ Purdue University, West Lafayette, Indiana - USA, 47907}
\begin{document}
\date{}
\maketitle

\begin{abstract}
This study presents a generalized multiscale multimesh finite element method (M${}^2$-FEM) that addresses several long-standing challenges in the numerical simulation of integral structural theories, often used to model multiscale and nonlocal effects. The major challenges in the numerical simulation of integral boundary value problems are primarily rooted in the coupling of the spatial discretization of the global (parent) and integral (child) domains which severely restricts the computational efficiency of existing algorithms by imposing an implicit trade-off in the accuracy achieved by the child domain and in the resources dedicated to the simulation of the overall parent domain. One of the most defining contributions of this study consists in the development of a mesh-decoupling technique that generates isolated sets of meshes such that the parent and child domains can be discretized and approximated independently. This mesh-decoupling has a multi-fold impact on the simulation of integral theories such that, when compared to existing state-of-the-art techniques, the proposed algorithm achieves simultaneously better numerical accuracy and efficiency (hence allowing a greater flexibility in both mesh size and computational cost trade-off decisions), greater ability to adopt generalized integral kernel functions, and the ability to handle non-regular (non-rectangular) domains via unstructured meshing. In this study, we choose a benchmark problem based on an extended version of the Eringen's nonlocal elasticity theory (implicitly, a multiscale theory) that leverages the use of generalized attenuation kernels and non-constant horizons of nonlocality. Nonetheless, the proposed $\text{M}^2$-FEM algorithm is very general and it can be applied to a variety of integral theories, even beyond structural elasticity.

\subsection*{Highlights}
\begin{itemize}[font=\small,leftmargin=*]
    \item Finite element solver for integral structural theories: $\text{M}^2$-FEM
    \item Mesh-decoupling and scale-bridging for scale-localized approximations 
    \item Ability to achieve both high accuracy and high efficiency simultaneously
    \item Generalized kernel functions allowing modeling different material classes
    \item Unstructured mesh to simulate irregular structures and non-similar integral horizons
\end{itemize}

\subsection*{Keywords: Integral elasticity | Multiscale effects | Nonlocal effects | Mesh-decoupling | Scale-bridging}
\end{abstract}

\clearpage
\section{Introduction}
\label{sec: introduction}
In recent years, several complex classes of heterogeneous material systems including, for example, metamaterials \cite{ganti2020topological}, porous solids \cite{mannan2017correlations}, functionally-graded solids \cite{parameswaran2000processing}, lattice structures \cite{russillo2022wave}, and composites \cite{kanoute2009multiscale}, have become progressively more accessible from a fabrication perspective. These novel classes of materials find important applications in the most diverse fields of engineering, nanotechnology, biotechnology, and medicine. The ability to exploit these complex material systems in practical applications requires structural models that are capable of reliably and accurately predicting the structural response. Several theoretical and experimental investigations have reinforced the understanding that these materials exhibit non-negligible non-classical behavior, that manifest at the continuum level, as multiscale effects \cite{kanoute2009multiscale,hoekstra2014multiscale,van2020roadmap}, nonlocal effects \cite{patnaik2022role,mei2022nonlocality,madeo2016first}, and even combined multiscale-nonlocal effects \cite{ding2022multiscale,ding2022multiscale2}. 

The inability of the classical local continuum theory to efficiently capture the non-classical multiscale and nonlocal effects motivated the development of non-classical continuum theories. From a broad mathematical perspective, all existing non-classical continuum theories can be broadly identified either as integral or gradient continuum theories. Detailed theoretical investigations and a diverse range of applications based over the past 50 years, have made it increasingly evident that the integral approach is more general when compared to the gradient approach \cite{batra159misuse,patnaik2022displacement,patnaik2020generalized,polizzotto2001nonlocal}. In the following, we provide a high-level discussion on nonlocal and multiscale effects to briefly clarify this latter aspect and to substantiate the choice, made later in this study, to simulate specific integral theories. 

From a physical perspective, both nonlocal and multiscale effects are characterized by a collective and multilevel transfer of information, from the finer structural (material and/or geometric) scales to the coarser (often, continuum) scale. The collective information leading to nonlocal effects is associated to multiple long-range interactions typically distributed in parallel within a specific finer structural scale \cite{polizzotto2001nonlocal,eringen1984theory}; in the case of multiscale effects, they correspond to a hierarchical distribution of multiple structural scales operating in parallel \cite{mcdowell2010perspective,liu2021multilayer}. The coexistence of hierarchical structural scales and localized long-range connections at a given scale, leads to multiscale nonlocality \cite{ding2022multiscale}. 

Integral approaches allow resolving the collective information by continuously (in a mathematical sense) accounting the contributions via model-specific convolution(s) over a set of interacting points localized either in a pre-defined nonlocal horizon \cite{patnaik2020generalized}, or a set of structural scales \cite{mcdowell2010perspective}, or both simultaneously \cite{ding2022multiscale}. Irrespective of the specific constitutive material model, integral approaches introduce a kernel function that characterizes the relative decay in the multilevel exchange of information, and a horizon of influence (as a set of material points in a localized scale, or a set of scales, or both) that supports the kernel function. Note that, from a mathematical standpoint, the multilevel exchange of information originating either from multiple points (as in, nonlocal effects) and/or multiple hierarchical scales (as in, multiscale effects), results in a similar structure of the enriched integral constitutive equations at the continuum scale. Consequently, in literature, nonlocal effects are often treated implicitly as multiscale effects (and the vice versa) \cite{mcdowell2010perspective,polizzotto2001nonlocal,eringen1984theory}. With this understanding, in the following, we have used ``multiscale effects" as a broad term that also includes nonlocal effects.

Gradient approaches gained popularity until the past few decades mostly for their easier numerical implementation, when compared to the significant complexities involved in integral approaches. However, in contrast to the seamless information resolution afforded by the integral approaches, gradient approaches use a limited set of higher-order gradients that capture the non-classical information from a limited number of higher-order (i.e. nearest, next-nearest, and so on) interactions or structural scales. In fact, gradient approaches are often derived from integral approaches by adding specialized assumptions on the functional variation and symmetry of the kernel function \cite{mcdowell2010perspective,polizzotto2001nonlocal}. These assumptions introduce major theoretical constraints in the modeling of heterogeneous and asymmetric microstructures \cite{batra159misuse}, and generalized attenuation-dispersion behavior \cite{patnaik2022displacement}. Consequently, gradient approaches have found rather narrow and, on several instances, even erroneous applications \cite{batra159misuse}. Thanks to the ability of admitting generalized kernels, integral approaches (particularly the more recent approaches) enable the modeling of a wide range of microstructural classes \cite{patnaik2022distillation,stempin2021formulation,failla2020advanced,nair2019nonlocal} exhibiting intricate localized (softening or stiffening) responses \cite{ding2022multiscale} and various attenuation-dispersion characteristics \cite{patnaik2020generalized,russillo2022wave}. Indeed, in recent times, integral approaches have received renewed and greater attention in non-classical modeling.

A detailed literature survey reveals that, despite the growing attention, numerical techniques necessary to fully realize the potential of integral frameworks, and complement their application, are yet to be developed \cite{van2020roadmap,chopard2014framework}. In \S\ref{sec: Motivation} below, we discuss the limitations of existing numerical techniques and their common origin in the sharing of mesh discretization among the different structural scales. Following this discussion, we present an overview of the main contributions of the present study to the development of advanced numerical techniques applicable to integral frameworks, in \S\ref{sec: Contributions}.

\clearpage
\subsection{Limitations of existing numerical solvers for integral theories}
\label{sec: Motivation}
To-date, only a handful of numerical techniques, including the NL-FEM (nonlocal finite element method) \cite{pisano2009nonlocal,sidhardh2018effect} and the f-FEM (fractional-order finite element method) \cite{patnaik2020ritz,patnaik2021fractional}, have successfully simulated generalized integral structural theories. The most critical limitation, common to these existing numerical strategies, consists in the explicit sharing of a globally defined spatial mesh in the discretization of the coarse (parent) and fine (child) length-scales (or domains). The shared discretization, or equivalently, mesh-coupling, between different material scales adversely impacts the simulation of multiscale problems leading to:
\begin{itemize}[leftmargin=*]
    \item \textit{Inability to deliver high accuracy and efficiency, simultaneously}. Mesh coupling increasingly leads to computationally very expensive simulations with minimal increase in accuracy. Note that the overall accuracy in multiscale problems is strongly dependent on the spatial resolution achieved at the child domain. Consequently, the finest spatial resolution directly determines the size of the overall mesh, as a result of the underlying spatial coupling. Hence, in applications involving a greater disparity between coarse and fine scales, this restriction could lead to a very computationally expensive mesh at the level of the coarse domain, potentially even preventing the ability to numerically solve the problem. As an example, consider a porous beam with very fine pores' scales compared to the characteristic length of the beam. In this case, it is not difficult to envision how a fine discretization for the child domain could lead to a multi-fold increase in the mesh size (in terms of the number of elements) at the level of the parent (i.e. beam scale) domain, drastically impacting the overall efficiency of the simulation technique.
    
    \item \textit{Inability to simulate irregular structural configurations}. The direct sharing of the spatial discretization implicitly requires the parent and child domains to be geometrically similar which, in practice, restricts the application of existing approaches to the simulation of multiscale structures with either rectangular domains (such as, for example, beams and plates) or at best, structures that can be reduced parametrically to rectangular domains (such as rectangular panels and shells). Consequently, it is not surprising that the application of existing approaches to multiscale structures with irregular material geometry and (geometrically) non-similar integral (wither multiscale or nonlocal) horizon geometries are absent from the literature.
\end{itemize}

Notably, in addition to the mesh-coupling, existing approaches are not well equipped to directly accommodate the most generalized (singular) kernel functions such as, for example, the power-law kernel (in constant- \cite{patnaik2020generalized}, variable- \cite{patnaik2021variable}, or distributed-order \cite{ding2022multiscale} version) which has found widespread application in the recent times \cite{patnaik2022distillation,stempin2021formulation,failla2020advanced}. f-FEM, that simulates differ-integral partial differential equations based on fractional calculus\footnote{Fractional calculus is a branch of mathematics that enables the differentiation and integration of functions beyond the classical set of natural numbers and extends it to real and complex numbers, and even functions. A discussion of fractional-order derivatives and integrals is beyond the scope of this study and the interested reader is referred to \cite{das2011functional} for detailed discussions.}, does admit singular kernels, but it still adopts coupled meshing.

The several limitations discussed above largely restrict the application of existing computational strategies to the simulation of highly multiscale structures characterized by strong underlying heterogeneity \cite{ostoja2002microstructural,mcdowell2010perspective}. An analysis of the most recent developments in the foundational theory and application of multiscale models reveals that classical phenomenological approaches (that can be more readily simulated by existing solvers) are of limited practical significance in the several novel classes of materials mentioned previously. The need for accuracy has suggested embracing the multiscale heterogeneous nature of these materials \cite{patnaik2022distillation,ding2022multiscale2,romkes2004adaptive}, as opposed to a more classical homogeneous and symmetric view typically associated with phenomenological (representative volume element, RVE-based) approaches. Some notable studies that provide concrete strategies for this renewed effort include mesoscale approaches informed by extensive molecular dynamics simulations \cite{liu2021multilayer,mcdowell2010perspective}, nonlocal approaches based on variable-order \cite{patnaik2021variable} and distributed-order \cite{ding2022multiscale} fractional calculus, and even, data-driven methods \cite{you2021data,flaschel2021unsupervised}. 

The large material parameter space associated with the aforementioned advanced theories ultimately necessitates the application of advanced inverse techniques (e.g. machine learning-based techniques) for constitutive discovery \cite{gallet2022structural,you2021data,patnaik2021variable,flaschel2021unsupervised}. The extensive datasets used in these discovery processes, are often generated artificially by the associated numerical solvers based either on \textit{a priori} knowledge of selective parameter combinations \cite{you2021data,jokar2021finite} or via novel strategies to implement one-shot discovery \cite{patnaik2022distillation,flaschel2021unsupervised}. To this end, the prohibitive computational cost required for better accuracy, the limited ability in modeling various geometries resulting from the mesh-coupling, and the inability in admitting generalized kernel functions, drastically impact the generation of critical training datasets and in general, even the forward implementation of these advanced theories in future applications. Evidently, existing multiscale computational approaches are not well-suited for the next-generation of applications. 

\subsection{Overall approach and contributions of the present study}
\label{sec: Contributions}
The previous discussion motivates the need to revisit the mesh topology and sharing strategy adopted in existing solvers. It is evident that the ideal computational approach, which addresses the limitations outlined in \S\ref{sec: Motivation}, requires an effective decoupling between the spatial discretization of the parent and child domains. It is important to understand that while the decoupling takes place in the numerical discretization, from a theoretical standpoint the information exchange across scales is still maintained. This condition can be achieved independently of the meshing strategy by ensuring that the adopted multiscale theory is thermodynamically consistent \cite{polizzotto2003unified,mcdowell2010perspective}, and by imposing appropriate connectivity among co-existing multiscale sets of material points in the numerical simulation (more commonly known as scale-bridging) \cite{chopard2014framework,chung2021nonlocal}. These concepts define the broad approach of the present study.

\textit{One of the most significant contribution of this study, towards advanced integral simulation, is the development of a novel mesh-decoupling technique that lays the foundation of a multiscale multimesh finite element method, denoted as $\text{M}^2\text{-FEM}$}. 

Recall that the weak-form approximation of integral theories consists of the approximation of (at least) two nested integrals, the parent-domain integral and the child-domain constitutive (multiscale/nonlocal) integral. The proposed mesh-decoupling technique generates isolated sets of meshes in the parent and child domains. By approximating the nested constitutive integral in the decoupled child mesh, the integral constitutive information is encapsulated within the child domain. This approach minimizes the transmission of numerical error generated from the child domain approximation to the parent domain (and likewise in the reverse direction). Thus, mesh-decoupling ultimately decouples the errors in the approximation of the two domain-specific integrals. The most significant outcome of this decoupling is the remarkable ability of $\text{M}^2$-FEM to independently tune parent and child mesh sizes, to achieve the desired efficiency and accuracy, simultaneously. Further, decoupling also allows unstructured meshing of both domains to model irregular configurations, and tailored numerical algorithms to embrace generalized singular kernels while approximating child domain integrals. In the present work, the development of $\text{M}^2\text{-FEM}$ is complemented by the following additional set of studies: 
\begin{itemize}[leftmargin=*]
    \item \textit{Validation and convergence analyses} are performed through a series of targeted case studies to validate the implementation of the mesh-decoupling mechanism and establish convergence conditions. The role of the child and parent domain resolutions, and the kernel function in determining the accuracy and convergence are investigated. As we will show, $\text{M}^2$-FEM enables a rapidly convergent and accurate simulation with relatively coarse parent domain mesh; a remarkable feature that is not attainable with existing solvers.
    
    \item \textit{Computational complexity analysis and simulation of non-rectangular structures} are presented to showcase the potential impact of $\text{M}^2$-FEM on the next-generation of multiscale problems. We will perform specific case studies to demonstrate how the superior convergence characteristics enable a relative tuning of the parent and child mesh resolutions to achieve superior accuracy and faster simulation times, simultaneously. Finally, a series of advanced geometries, identified from practically relevant structures in a variety of fields, are considered to showcase the abilities of $\text{M}^2$-FEM; existing solvers cannot simulate the considered geometries.
\end{itemize}

For all the simulations presented in this study, we have adopted as the underlying structural theory an extended version of Eringen's integral (strain-driven) nonlocal theory \cite{eringen1984theory} that admits non-singular and singular kernel functions. Note that Eringen's nonlocal theory is implicitly multiscale in nature \cite{mcdowell2010perspective}. To this end, we emphasize that the algorithm developed is very general and can be directly applied to any integral multiscale formulation related to structural mechanics or other continuum theories. We decided to adopt Eringen's approach because the broader scientific community still relates more closely to this seminal theory. It is also important to note that critical aspects, such as model convexity, thermodynamic consistency, and mathematical well-posedness are difficult to address when using the strain-driven approach \cite{polizzotto2001nonlocal}. We merely note that the recently proposed displacement-driven approach \cite{patnaik2022displacement} is a more suitable candidate since it is intrinsically convex, thermodynamically consistent, and well-posed. Nonetheless, since model convexity is essential for finite element solvers, the successful simulation of the strain-driven theory in this study, essentially guarantees a rather straightforward application of $\text{M}^2$-FEM also to the simulation of the displacement-driven theory.

The remainder of this paper is structured as follows. The mesh-decoupling technique and the scale-bridging algorithms are developed in \S\ref{sec: 2}, alongside the finite element approximations of both parent and child domains. Next, in \S\ref{sec: 3}, we develop the matrix representation of the approximated components, assemble the domain-specific approximations, and derive the algebraic expressions for computer implementation. Results of validation and convergence studies are presented in \S\ref{sec: 4}. Simulations obtained from a series of practically relevant multiscale configurations are presented in \S\ref{sec: 5} to highlight the envisioned impact of $\text{M}^2$-FEM. Finally, a brief summary and the future outlook are presented in \S\ref{sec: 6}.

\section{Mathematical formulation: decoupled multiscale discretization and numerical approximation}
\label{sec: 2}
In this section, we develop the mesh-decoupling algorithm that is at the foundation of M${}^2$-FEM. This development includes the definition of separate meshes for the discretization of different material scales as well as the use of the different meshes to derive the multiscale finite element approximations for the simulation of the background theory. Note that, while the algorithm is general, we will derive the approximations for a 2D problem (as opposed to the full 3D problem) to facilitate the understanding of important aspects of the algorithm. Here below, we briefly summarize Eringen's integral theory (focusing on the constitutive relations and governing equations) and the notation scheme adopted in the study, to more efficiently utilize them in deriving and discussing the discretization and approximation techniques.

\subsection{Brief summary of the adopted integral theory}
\label{sec: Brief_theory}
According to the strain-driven approach to nonlocal elasticity\cite{eringen1984theory}, the stress-strain constitutive relation at any material point in a nonlocal solid, can be expressed using integral operators in the following form:
\begin{equation}\label{eq: constitutive_relation}
    \bm{\sigma} = \bm{C}:\underbrace{\int_{\mathcal{H}} \mathcal{K}(\bm{x},\bm{x}^\prime) \bm{\epsilon}(\bm{x}^\prime) \mathrm{d}\bm{x}^\prime}_{\textrm{nonlocal integral}}
\end{equation}
where $\bm{\sigma}$, $\bm{C}$, and $\bm{\epsilon}$ are the nonlocal stress, stiffness, and local strain tensors, respectively. Further, $\bm{x}$ denotes a specific material point on the nonlocal solid (say, $\Omega$); and $\bm{x}^\prime \in \mathcal{H}$ is a dummy spatial variable denoting the (remaining) materials points that interact with $\bm{x}$ via long-range connections and are contained within a characteristic nonlocal horizon $\mathcal{H}$. The strength of the nonlocal effects are captured by the nonlocal kernel function denoted as $\mathcal{K}(\bm{x},\bm{x}^\prime)$ with support on the nonlocal horizon $\mathcal{H}$.

The strong form of the governing equations can be derived by variational minimization of the potential energy functional obtained from the constitutive framework presented in Eq.~(\ref{eq: constitutive_relation}), as \cite{polizzotto2001nonlocal}: 
\begin{equation}\label{eq: strong_governing_eqn}
\begin{aligned}
    \bm{\nabla}\cdot\bm{\sigma}+\rho\bm{f}-\rho\bm{\Ddot{u}}&=0,\quad \forall \bm{x} \in \bm{\Omega}\\
\end{aligned}
\end{equation}
The governing equations are subject to the following boundary and initial conditions:
\begin{equation}\label{eq: BCs}
\begin{aligned}
    \delta\bm{u}=0 \quad \text{or} \quad \bm{\sigma}\cdot\bm{n}-\bm{t}=0, \quad \forall \bm{x} \in \bm{\Gamma}\\
    \delta\bm{u}=0 \quad \text{and} \quad \delta\dot{\bm{u}}=0, \quad \forall \bm{x}(t=0) \\
\end{aligned}
\end{equation}
where $\bm{\Gamma}$ denotes the boundary of the nonlocal solid, and $\bm{t}$ denotes the traction force acting on the surface whose outward normal is indicated by $\bm{n}$. Analogous to the classical approach, we use the weak form representation of Eq.~(\ref{eq: strong_governing_eqn}), presented below (see \cite{polizzotto2001nonlocal} for derivation), for finite element discretization and approximation:
\begin{equation}\label{eq: weak_governing_eqn}
\begin{aligned}
    \int_{\bm{\Omega}} \bm{\sigma}:\bm{\nabla}\delta\bm{u} \mathrm{d}\bm{\Omega} 
    - \int_{\bm{\Omega}} \rho(\bm{f}-\Ddot{\bm{u}})\delta \bm{u} \mathrm{d}\bm{\Omega} 
    - \int_{\partial\bm{\Omega}} \bm{\sigma}\cdot\bm{n}\delta \bm{u} \mathrm{d}\bm{\Gamma} &= 0 \\
\end{aligned}
\end{equation}
In proceeding further, we assume that the nonlocal kernel function $\mathcal{K}(\bm{x},\bm{x}^\prime)$ is positive-definite and symmetric such that the total strain energy functional is convex and the governing equations in Eq.~(\ref{eq: strong_governing_eqn}) are well-posed \cite{polizzotto2001nonlocal}.

It is evident from the Eqs.~(\ref{eq: constitutive_relation},\ref{eq: weak_governing_eqn}) that the finite element implementation of the assumed background theory requires the approximations of different physical variables defined on separate material scales. These approximations correspond either to: 1) the weak-form integral over the material domain ($\bm{\Omega}$), and 2) the convolution integral defined over the nonlocal horizon ($\mathcal{H}$). Henceforth, we simply refer to the two material scales as the parent domain (that is, the material domain) and the child domain (that is, the nonlocal horizon). 

\subsection{Remarks on the adopted notation scheme}
\label{sec: 2.2}
The separate meshing schemes used to discretize the two different domains and the integral nature of the underlying formulation result in a large number of variables and finite element parameters. In the interest of streamlining the following derivations and discussions, as well as to facilitate the reading, we discuss the notation scheme here below. The most general notation for a given variable, say $\square$, that encompasses the subsequent development can be expressed as:
\begin{equation}\label{eq: notation_form}
    \prescript{\mathrm{u_-}}{l_-}{}\square^{\mathrm{u_+}}_{l_+}
\end{equation}
where the sets of superscripts $\{\mathrm{u}_-,\mathrm{u}_+\}$ and subscripts $\{l_-,l_+\}$ can be specified to denote different variables. 

We typeset the superscripts in Roman font and use them to generate symbols for different variables well-defined either on the decoupled-meshes or on the solid domain. The complete list of superscripts is here below:
\begin{itemize}[leftmargin=*]
    \item $\{\mathrm{p},\mathrm{c}\}$ represent quantities defined on the parent and child domain.
    \item $\{\mathrm{b},\mathrm{u},\mathrm{t}\}$ represent quantities defined on the complete boundary, and the segments of the boundary subjected to displacement and traction loads.
\end{itemize}
For enhanced readability, we typeset the subscript using italic fonts and use them for indexing. Besides the commonly used indices $i$ and $j$, we define the following indices to represent specific physical quantities:
\begin{itemize}[leftmargin=*]
    \item $\{r,\tilde{r},s\}$ are used for parent, dummy parent (defined in $\mathcal{H}$, see \S\ref{sec: 2.4.2}), and child element indexing.
    \item $\{d,g\}$ are used for nodal and Gauss point indexing.
    \item $\{u,v\}$ are used for number of degrees of freedom (DOF) and weighting function indexing.
\end{itemize}
Further, we use the right-hand scripts $\{\mathrm{u}_+,l_+\}$ to indicate quantities defined at the current scale (for example, $\Omega_r^\textrm{p}$ and $\Omega_s^\textrm{c}$ that define finite elements in the parent domain and child domain, respectively) and left-hand scripts $\{\mathrm{u}_-,l_-\}$ to indicate coarser scale quantities. This convention proves effective when describing multiscale quantities. Consider for example, $\prescript{\textrm{p}}{rg}{}\Omega_{s}^\textrm{c}$, which represents finite elements in the child domain (see \S\ref{sec: 2.4.1} or Eq.~(\ref{eq: child_meshing})). In $\prescript{\textrm{p}}{rg}{}\Omega_{s}^\textrm{c}$, the right scripts ${\square}_{s}^\textrm{c}$ specify that $\prescript{\textrm{p}}{rg}{}\Omega_{s}^\textrm{c}$ is the $s^\mathrm{th}$ element in the child domain, and the left scripts $\prescript{\textrm{p}}{rg}{\square}$ specify that this child element is defined at the nonlocal horizon of $g^\mathrm{th}$ parent Gauss point in the $r^\mathrm{th}$ parent element.

Two final remarks on the notation: 1) we have used a comma to denote partial derivatives in agreement with the classical tensor algebra notation \cite{oden2012introduction}; and 2) in the interest of a general derivation, we have used different symbols, $\Phi(\cdot)$ and $\Psi(\cdot)$, to denote the shape functions used to interpolate the geometrical coordinates and the displacement degrees of freedom, respectively; the remaining variables are introduced as they appear in the derivation below.

\subsection{Parent domain discretization and approximation}
\label{sec: 2.3}
We first develop the finite element approximation of the weak-form governing equations over the parent domain $\bm{\Omega}$ and boundary $\bm{\Gamma}$. Note that, since the governing equations on the material scale presented in Eq.~(\ref{eq: strong_governing_eqn}) exhibit a functional character similar to their classical (local) counterpart, the parent domain approximations very closely follows a classical finite element approximation of the 2D plane strain problem. However, we fully derive the parent domain approximations, following a discretization of the weak-form integral at the parent level in order to: 1) clearly highlight the mesh-decoupling in action and the connectivity algorithm developed to preserve the multiscale exchange of information, and 2) enable an efficient discussion on the matrix representation and integration of the dual-level approximations in the final implementation in \S\ref{sec: 3}.

\subsubsection{Discretization}
\label{sec: 2.3.1}
Consider the nonlocal elastic material occupying, in its undeformed state, a 2D Euclidean parent domain $\bm{\Omega}^\textrm{p}$ referred to orthogonal Cartesian coordinates $\bm{x}=(x,y)$. The 2D domain is discretized into $N^\textrm{p}$ disjoint finite elements $\Omega_r^\textrm{p}$ such that:
\begin{equation}
    \bigcup_{r=1}^{N^\textrm{p}}\Omega_r = \bm{\Omega}, \quad \Omega_i \cap \Omega_j = \varnothing, ~\forall i \neq j
\end{equation}
where the superscript $\textrm{p}$ indicates quantities related to the $\textit{parent domain}$. In all the computational developments considered in the following, the parent mesh is assumed to be non-rectangular in shape, unless specified otherwise. The boundary $\bm{\Gamma}$ is also discretized into $N^{\textrm{b}}$ elements that can be differentiated into the following different types depending on the nature of the boundary condition:
\begin{equation}
\begin{aligned}
    \bigcup_{r=1}^{N^\textrm{b}}\Gamma_r &= \bigcup_{i=1}^{N^\textrm{u}}\Gamma_i^\textrm{u} + \bigcup_{j=1}^{N^\textrm{t}}\Gamma_j^\textrm{t} = \bm{\Gamma}, \quad &&\Gamma_i^\textrm{u} \cap \Gamma_j^\textrm{t} = \varnothing, ~\forall i \in [1,N^\textrm{u}],~j \in [1,N^\textrm{t}] \\
    \bigcup_{r=1}^{N^\textrm{u}}\Gamma_r^\textrm{u} &= \bm{\Gamma}^\textrm{u}, \quad &&\Gamma_i^\textrm{u} \cap \Gamma_j^\textrm{u} = \varnothing, ~\forall i \neq j \\
    \bigcup_{r=1}^{N^\textrm{t}}\Gamma_r^\textrm{t} &= \bm{\Gamma}^\textrm{t}, \quad &&\Gamma_i^\textrm{t} \cap \Gamma_j^\textrm{t} = \varnothing, ~\forall i \neq j \\
\end{aligned}
\end{equation}
where $N^\textrm{u}$ and $N^\textrm{t}$ are the number of boundaries subjected to the displacement boundary condition (represented by $\bm{\Gamma}^\textrm{u}$) and traction boundary conditions (represented by $\bm{\Gamma}^\textrm{t}$), respectively.

Using the above parent domain discretization, the continuous integral corresponding to the simplified variational governing equations (see Eq.~(\ref{eq: weak_governing_eqn})) is approximated as the limit of the following summation:
\begin{equation}
    \sum_{r=1}^{N^\textrm{p}}\int_{\Omega_r^\textrm{p}} \bm{\sigma} : \bm{\nabla}_{\bm{x}}\delta\bm{v} \mathrm{d}\bm{\Omega}
    - \sum_{r=1}^{N^\textrm{p}}\int_{\Omega_r^\textrm{p}} \rho(\bm{f}-\Ddot{\bm{u}})\cdot\delta\bm{v} \mathrm{d}\bm{\Omega}
    - \sum_{r=1}^{N^\textrm{t}}\int_{\Gamma_r^\textrm{t}} \bm{\sigma}\cdot\bm{n} \delta\bm{v} \mathrm{d}\bm{\Gamma} = 0
\end{equation}
where $\bm{u}$ and $\delta\bm{v}$ denote the nodal degrees of freedom (DOF) and the weighting function. Further, the boundary conditions are expressed as:
\begin{equation}
\begin{aligned}
    \bm{u}|_{\bm{\Gamma}^\textrm{u}} &= \overline{\bm{u}} \\
    \bm{\sigma}\cdot\bm{n}|_{\bm{\Gamma}^\textrm{t}} &= \overline{\bm{t}} \\
\end{aligned}
\end{equation}
where $\overline{\bm{u}}$ and $\overline{\bm{t}}$ denote the externally applied displacement and traction loads. By combining the approximations presented above for the governing equations and boundary conditions, the first step discretization developed here at the parent domain $\bm{\Omega}^\textrm{p}$ results in the following approximation of the governing variational equation:
\begin{equation}\label{eq: governing_eqn_parent_domain}
\begin{aligned}
    \sum_{r=1}^{N^\textrm{p}}\int_{\Omega_r^\textrm{p}} \bm{\sigma} : \bm{\nabla}_{\bm{x}}\delta\bm{v} \mathrm{d}\Omega_r^\textrm{p}
    - \sum_{r=1}^{N^\textrm{p}}\int_{\Omega_r^\textrm{p}}\rho(\bm{f}-\Ddot{\bm{u}})\cdot\delta\bm{v} \mathrm{d}\Omega_r^\textrm{p}
    - \sum_{r=1}^{N^\textrm{t}}\int_{\Gamma_r^\textrm{t}} \overline{\bm{t}}\cdot\delta\bm{v} \mathrm{d}\Gamma_r^\textrm{t} &= 0
\end{aligned}
\end{equation}

\subsubsection{Numerical integration}
\label{sec: 2.3.2}
As evident from Eq.~(\ref{eq: governing_eqn_parent_domain}), the approximation of the variational equation consists of two different components: 1) the numerical evaluation of the discretized integral expression over the parent domain elements $\Omega_r^\textrm{p}$, and 2) the numerical evaluation of the discretized integral expression over the boundary elements $\Gamma_r^\textrm{t}$. In the following, we first evaluate the parent domain integrals and then consider the boundary domain integrals.

In order to perform the numerical integration, we adopt an isoparametric formulation and introduce a 2D natural coordinate system. In more specific terms, for a given parent element $\Omega_r^\textrm{p}$, $\forall r \in [1,N^\textrm{p}]$, we transform the material points $\bm{x}=(x,y)$ in $\Omega_r^\textrm{p}$ from the 2D Cartesian coordinates to natural coordinates $\bm{\xi}_r^\textrm{p}=(\xi_r^\textrm{p},\eta_r^\textrm{p})$. This isoparametric mapping can be expressed as:
\begin{equation}\label{eq: affine_transform_parent_element}
    \bm{x}(\bm{\xi}_r^{\textrm{p}}) = \sum_{d=1}^{d^{\textrm{p}}} \Phi_{d}^\textrm{p}(\bm{\xi}_r^\textrm{p})\bm{x}_{rd}^\textrm{p}
    = \{\Phi_d(\bm{\xi}_r^\textrm{p})\}^{\textrm{T}} \{\bm{x}_{rd}^\textrm{p}\}
\end{equation}
such that the parent element $\Omega_r^\textrm{p}$ can be transformed into a rectangular element $\hat{\Omega}_r^\textrm{p}$ in the natural coordinate system. In the latter transformation, $\Phi_{d}^\textrm{p}(\cdot)$ and $d^\textrm{p}$ are the shape functions and the number of nodes in the parent element $\Omega_r^\textrm{p}$, respectively. Finally, $\{\Phi_d(\bm{\xi}_r^\textrm{p})\}$ and $\{\bm{x}_{rd}^\textrm{p}\}$ (where $d=1,2,...,d^\textrm{p}$), that denote the vectors assembling the nodal shape functions and the corresponding nodal coordinates, respectively, are defined as:
\begin{subequations}
\begin{align}
    \{\Phi_d(\bm{\xi}_r^\textrm{p})\} &= \{\Phi_1(\bm{\xi}_r^\textrm{p}),\Phi_2(\bm{\xi}_r^\textrm{p}),...,\Phi_{d^\textrm{p}}(\bm{\xi}_r^\textrm{p})\}^\textrm{T} \\
    \{\bm{x}_{rd}^\textrm{p}\} &= \left\{\{x_{r}^\textrm{p},y_{r}^\textrm{p}\}_1^\textrm{T},\{x_{r}^\textrm{p},y_{r}^\textrm{p}\}_2^\textrm{T},...,\{x_{r}^\textrm{p},y_{r}^\textrm{p}\}_{d^\textrm{p}}^\textrm{T}\right\}^\textrm{T}
\end{align}
\end{subequations}
where the subscript $r$ denote quantities supported on the $r^\textrm{th}$ parent element. Note that, analogous to classical finite element methodology, the isoparametric transformation between the Cartesian and the natural coordinate systems is affine in nature \cite{oden2012introduction}. More specifically, the shape functions $\Phi_d^\textrm{p}(\cdot)$ used for geometrical interpolation should be linear functions when defined in either coordinate systems; see additional details in \S\ref{sec: 2.4.2}. We consider Gauss points $\bm{\xi}_{rg}^\textrm{p}=(\xi_{rg}^\textrm{p},\eta_{rg}^\textrm{p})$ and the corresponding weights $w_{rg}^\textrm{p}$ in the isoparametric element $\hat{\Omega}_r^\textrm{p}$. Now by using the classical Gauss-Legendre quadrature principles, the discretized integrals in Eq.~(\ref{eq: governing_eqn_parent_domain}) over the parent element $\Omega_r^\textrm{p}$, with Gauss points $\bm{\xi}_{rg}^\textrm{p}=(\xi_{rg}^\textrm{p},\eta_{rg}^\textrm{p})$ and the corresponding weights $w_{rg}^\textrm{p}$, can be approximated as:
\begin{subequations}\label{eq: Gauss_integration_parent_element}
\begin{align}
    \int_{\Omega_r^\textrm{p}} \bm{\sigma} : \bm{\nabla}_{\bm{x}}\delta\bm{v} \mathrm{d}\Omega_r^\textrm{p} 
    &= \sum_{g=1}^{g^\textrm{p}} {w}_{rg}^\textrm{p} \vert\bm{J}_{r}^\textrm{p}(\bm{\xi}_{rg}^\textrm{p})\vert \bm{\sigma}(\bm{\xi}_{rg}^\textrm{p}) : \bm{\nabla}_{\bm{x}}\delta\bm{v}(\bm{\xi}_{rg}^\textrm{p}) \\
    \int_{\Omega_r^\textrm{p}}\rho(\bm{f}-\Ddot{\bm{u}}) \delta\bm{v} \mathrm{d}\Omega_r^\textrm{p}  
    &= \sum_{g=1}^{g^\textrm{p}} w_{rg}^\textrm{p} \vert\bm{J}_{r}^\textrm{p}(\bm{\xi}_{rg}^\textrm{p})\vert \rho(\bm{\xi}_{rg}^\textrm{p})(\bm{f}(\bm{\xi}_{rg}^\textrm{p})-\Ddot{\bm{u}}(\bm{\xi}_{rg}^\textrm{p})) \cdot \delta\bm{v}(\bm{\xi}_{rg}^\textrm{p})
\end{align}
\end{subequations}
where $g^\textrm{p}$ denotes the number of Gauss points. $\vert\bm{J}_{g}^\textrm{p}(\bm{\xi}_{rg}^\textrm{p})\vert$ denotes the Jacobian of the transformation at the Gauss point $\bm{\xi}_{rg}^\textrm{p}$ and is evaluated in the following manner:
\begin{equation}
    \bm{J}_{r}^\textrm{p}(\bm{\xi}_{rg}^\textrm{p})
    = \frac{\partial \bm{x}(\bm{\xi}_{rg}^\textrm{p})}{\partial \bm{\xi}_{r}^\textrm{p}} 
    = \sum_{d=1}^{d^\textrm{p}} \frac{\partial \Phi_{d}^\textrm{p}(\bm{\xi}_{rg}^\textrm{p})}{\partial \bm{\xi}_{r}^\textrm{p}} \bm{x}_{rd}^\textrm{p}
    = \sum_{d=1}^{d^\textrm{p}} \bm{x}_{rd}^\textrm{p} \otimes \Phi_{d,\bm{\xi}_r^\textrm{p}}^\textrm{p}(\bm{\xi}_{rg}^\textrm{p})
    = \{\Phi_{d,\bm{\xi}_r^\textrm{p}}(\bm{\xi}_{rg}^\textrm{p})\}^\textrm{T} \{\bm{x}_{rd}^\textrm{p}\}
\end{equation}
in which $\otimes$ denotes the outer product between two tensors (note that the subscripts are indices of $\bm{x}$ and $\Phi$, not indices of components of these tensors). As we have clarified in \S\ref{sec: 2.2}, the subscript $\bm{\xi}_r^\textrm{p}$ (that consists of $\xi_r^\textrm{p}$ and $\eta_r^\textrm{p}$) indicates the partial derivative with respective to the isoparametric variable in the $r^\textrm{th}$ parent domain.

Note that the expressions derived in Eq.~(\ref{eq: Gauss_integration_parent_element}) require an evaluation of the displacement ($\Ddot{\bm{u}}$) and the weighting function ($\delta\bm{v}$ and $\bm{\nabla}_{\bm{x}}\delta\bm{v}$). By adopting the classical Galerkin approach, the weighting functions are chosen as the shape functions adopted to interpolate the displacement variables. Consequently, the displacement and weighting functions can be interpolated in terms of their corresponding nodal quantities in the following manner:
\begin{subequations}\label{eq: u_and_delta_v}
\begin{align}
    \Ddot{\bm{u}}(\bm{\xi}_{rg}^\textrm{p}) &= \sum_{u=1}^{m^\textrm{p}} \Psi_{u}^\textrm{p}(\bm{\xi}_r^\textrm{p})\Ddot{\bm{u}}_{ru}^\textrm{p} \\
    \delta\bm{v}(\bm{\xi}_{rg}^\textrm{p}) &= \sum_{v=1}^{m^{\textrm{p}}} \Psi_{v}^p(\bm{\xi}_r^p)\delta\bm{v}_{rv}^\textrm{p}
\end{align}
\end{subequations}
where subscripts $u,v \in \{1,2,...,m^\textrm{p}\}$ are the two indices for $\bm{u}$ and $\delta\bm{v}$, respectively; and $m^\textrm{p}$ is the total number of DOF within the parent element. 
Now, by using the shape function $\Psi(\cdot)$, the gradient of the weighting functions $\bm{\nabla}_{\bm{x}}\delta\bm{v}$ can be immediately approximated as:
\begin{equation}\label{eq: grad_delta_v}
    \bm{\nabla}_{\bm{x}} \delta\bm{v}(\bm{\xi}_{rg}^\textrm{p})
    = \sum_{v=1}^{m^{\textrm{p}}} \delta\bm{v}_{rv} \otimes \bm{\nabla}_{\bm{\xi}_{r}^\textrm{p}} \Psi_{v}^\textrm{p}(\bm{\xi}_{rg}^\textrm{p}) \frac{\partial\bm{\xi}_{rg}^\textrm{p}}{\partial\bm{x}}
    = \sum_{v=1}^{d^{\textrm{p}}} \delta\bm{v}_{rv} \otimes
     \left[\bm{\nabla}_{\bm{\xi}_{r}^\textrm{p}} 
    \Psi_{v}^\textrm{p}(\bm{\xi}_{rg}^\textrm{p}) (\bm{J}_{rg}^\textrm{p})^{-1}\right]
\end{equation}
where the italic letters $u$ and $v$ in the subscripts serve as indices to distinguish between shape functions for solutions (such as the term $\Ddot{\bm{u}}_{ru}^\textrm{p}$) and weight function terms (such as the term $\delta\bm{v}_{rv}^\textrm{p}$), respectively. 
Substituting Eqs.~(\ref{eq: u_and_delta_v},\ref{eq: grad_delta_v}) within Eq.~(\ref{eq: Gauss_integration_parent_element}), we finally obtain the approximation of the integral over $\Omega_r^\textrm{p}$ as:
\begin{subequations}\label{eq: Gauss_integration_parent_element_1}
\begin{align}
    \int_{\Omega_r^\textrm{p}} \bm{\sigma} : \bm{\nabla}_{\bm{x}}\delta\bm{v} \mathrm{d}\Omega_r^\textrm{p}
    &= \sum_{g=1}^{g^\textrm{p}} {w}_{rg}^\textrm{p} \vert\bm{J}_r^\textrm{p}(\bm{\xi}_{rg}^\textrm{p})\vert 
    \bm{\sigma}(\bm{\xi}_{rg}^\textrm{p}) : 
    \left[\sum_{v=1}^{d^{\textrm{p}}} \delta\bm{v}_{rv} \otimes \left[\bm{\nabla}_{\bm{\xi}_{r}^\textrm{p}} 
    \Psi_{v}^\textrm{p}(\bm{\xi}_{rg}^\textrm{p}) (\bm{J}_{rg}^\textrm{p})^{-1}\right]\right] \\
    \int_{\Omega_r^\textrm{p}}\rho(\bm{f}-\Ddot{\bm{u}}) \delta\bm{v} \mathrm{d}\Omega_r^\textrm{p} 
    &= \sum_{g=1}^{g^\textrm{p}} w_{rg}^\textrm{p} \vert\bm{J}_r^\textrm{p}(\bm{\xi}_{rg}^\textrm{p})\vert
    \rho(\bm{\xi}_{rg}^\textrm{p}) \left[\bm{f}(\bm{\xi}_{rg}^\textrm{p}) - \sum_{u=1}^{d^{\textrm{p}}} \Psi_{u}^\textrm{p}(\bm{\xi}_r^\textrm{p})\Ddot{\bm{u}}_{ru}^\textrm{p}\right]
    \cdot
    \left[\sum_{v=1}^{d^{\textrm{p}}} \Psi_{v}^\textrm{p}(\bm{\xi}_r^\textrm{p})\delta\bm{v}_{rv}^\textrm{p}\right]
\end{align}
\end{subequations}

By retracing the procedure outlined above, the integrals over the boundary elements in Eq.~(\ref{eq: governing_eqn_parent_domain}) are evaluated numerically, here below. For this purpose, the line elements in the Cartesian system (of the boundary) are mapped on to a 1D natural coordinate system via the following affine transformation:
\begin{equation}
    \bm{x}(\bm{\xi}_r^\textrm{t}) = \sum_{d=1}^{d^\textrm{t}} \Phi_d^\textrm{t}(\bm{\xi}_r^\textrm{t})\bm{x}_{rd}^\textrm{t}
\end{equation}
where $\Phi_d^\textrm{t}(\cdot)$ and $d^\textrm{t}$ are the $C^{0}$ Lagrangian shape function and the number of nodes in the boundary element $\Omega_r^\textrm{t}$, respectively. Again, by using Gauss-Legendre quadrature in isoparametric coordinates, the integral over the parent element $\Gamma_r^\textrm{t}$ shown in Eq.~(\ref{eq: governing_eqn_parent_domain}) can be approximated as:
\begin{equation}\label{eq: Gauss_integration_parent_boun}
    \int_{\Gamma_r^\textrm{t}} \overline{\bm{t}}\delta\bm{v} \mathrm{d}\Gamma_r^t  
    = \sum_{g=1}^{g^\textrm{t}} w_{rg}^\textrm{t}\vert\bm{J}_r^\textrm{t}\vert \overline{\bm{t}}(\bm{\xi}_{rg}^\textrm{t}) \delta\bm{v}(\bm{\xi}_{rg}^\textrm{t})
\end{equation}
where $g^\textrm{t}$ denotes the number of Gauss points at $\bm{\xi}_{rg}^\textrm{t}$ and the corresponding Gaussian weights $w_{rg}^\textrm{t}$. The Jacobian of the transformation, denoted above as $\vert\bm{J}_r^\textrm{t}\vert$, is obtained as:
\begin{equation}
    \bm{J}_r^\textrm{t} = \frac{\partial \bm{x}}{\partial \bm{\xi}_r^\textrm{t}} 
    = \sum_{d=1}^{d^\textrm{t}} \frac{\partial \Phi_{d}^\textrm{t}(\bm{\xi}_r^\textrm{t})}{\partial \bm{\xi}_r^\textrm{t}} \bm{x}_{rd}^\textrm{t}
    = \sum_{d=1}^{d^\textrm{t}} \bm{x}_{rd}^\textrm{t} \otimes \Phi_{d,\bm{\xi}_r^\textrm{t}}^\textrm{t}(\bm{\xi}_r^\textrm{t})
\end{equation}
Finally, by using the approximation of the weighting function term $\delta\bm{v}$ from Eq.~(\ref{eq: u_and_delta_v}), the boundary integrals over $\Gamma_r^\textrm{t}$ are evaluated as:
\begin{equation}\label{eq: Gauss_integration_parent_boun_1}
    \int_{\Gamma_r^\textrm{t}} \overline{\bm{t}}\delta\bm{v} \mathrm{d}\Gamma_r^\textrm{t}  
    = \sum_{g=1}^{g^\textrm{t}} w_{rg}^\textrm{t}\vert\bm{J}_r^\textrm{t}\vert \overline{\bm{t}}(\bm{\xi}_{rg}^\textrm{t}) 
    \left[\sum_{v=1}^{m^{\textrm{t}}} \Psi_{v}^\textrm{t}(\bm{\xi}_r^\textrm{t})\delta\bm{v}_{rv}^\textrm{t}\right]
\end{equation}
This concludes the finite element approximation at the parent domain level. 

\paragraph*{Resolution from parent to child domain.}\mbox{} Note that the evaluation of the stress $\bm{\sigma}$ requires the finite element approximation of the integral constitutive relation in Eq.~(\ref{eq: constitutive_relation}). The latter convolution integral requires an additional level of approximation that corresponds to the nonlocal horizon, that is, the child domain. This understanding brings us to the next step, that is the finite element discretization and approximations within the child domain.

\subsection{Child domain discretization and approximation}
\label{sec: 2.4}
In this section, we focus on the evaluation of nonlocal stress $\bm{\sigma}$ at each Gauss point in parent elements $\Omega_{rg}^\textrm{p}$ and boundary elements $\Gamma_{r}^\textrm{p}$. As evident from the integral constitutive relation in Eq.~(\ref{eq: constitutive_relation}), the numerical evaluation of the stress requires the information on the strain level at material points located within the nonlocal horizon $\mathcal{H}$. The evaluation of this nested integral poses several challenges when adopting an independent (child) mesh, decoupled from the parent mesh. Some of these challenges include the approximation of generalized (particularly, non-singular) nonlocal kernel functions and enforcing well-defined multiscale connectivity with the parent domain.

Before developing the discretization for the child mesh, we emphasize that the discretization and approximation of the child domain are directly influenced by the degree of singularity of the nonlocal kernel.  Non-singular kernels (for example, exponential kernels, typical of classical nonlocal approaches) lead to proper convolution integrals over the nonlocal horizon $\mathcal{H}$ that can be accurately approximated using the Gauss-Legendre quadrature principle. On the other hand, weakly-singular kernels (for example, power law kernels with exponents $\alpha\in(0,1)$, typical of fractional calculus based methods) lead to improper integrals that require the Gauss-Jacobi quadrature principle along with restrictions on the location of the quadrature points \cite{pang2013gauss}. Similarly, more advanced quadrature schemes could be used for strongly-singular nonlocal kernels (for example, power law kernels with exponents $\alpha>1$), these kernels are yet to find meaningful applications except in applications concerned with plastic damage such as, cracks and dislocations \cite{mcdowell2010perspective,liu2021multilayer}. In this study, we only consider elastic deformation across the scales and restrict the analysis to non-singular and weakly-singular nonlocal kernels.

\subsubsection{Discretization}
\label{sec: 2.4.1}
We consider the $r^{\mathrm{th}}$ parent element $\Omega_{r}^\textrm{p}$, the discretization of nonlocal horizon for each Gauss point $\bm{x}(\bm{\xi}_{rg}^\textrm{p})$, $g={1,2,...,g^\mathrm{p}}$ in $\Omega_{r}^\textrm{p}$, consists of the following two steps in the outlined sequence:
\begin{itemize}[leftmargin=*]
    \item \textit{Step 1:} determine the nonlocal horizon, that is the child domain $\prescript{\mathrm{p}}{rg}{}\bm{\Omega}^\mathrm{c}=\mathcal{H}(\bm{x}(\bm{\xi}_{rg}^\mathrm{p}))$, at the Gauss point. The prescript $\prescript{\textrm{p}}{rg}{}(\cdot)$ indicates quantities that belong to the $g^{\mathrm{th}}$ Gauss point in the $r^{\mathrm{th}}$ parent element.
    
    \item \textit{Step 2:} discretize the nonlocal horizon to obtain the child mesh.
\end{itemize}
Step 1 above determines the shape and geometry of the child domain $\prescript{\textrm{p}}{rg}{}\bm{\Omega}^\textrm{c}$ in the Euclidean space. Since the nonlocal horizon is contained within the material boundary $\bm{\Gamma}$, $\prescript{\textrm{p}}{rg}{}\bm{\Omega}^\textrm{c}$ at different Gauss points (particularly for Gauss points $\bm{x}(\bm{\xi}_{rg}^\textrm{p})$ near the material boundary) are truncated. Figure~(\ref{fig: unstructured_mesh}-b), which shows the mesh configuration of examples (which are simulated in \S\ref{sec: 4.1}), illustrates how the material boundary affects the shape of nonlocal horizon. 

Similar to parent domain, the 2D child domain $\prescript{\textrm{p}}{rg}{}\bm{\Omega}^\textrm{c}$ is discretized into $\prescript{\textrm{p}}{rg}{}N^\textrm{c}$ disjoint finite elements $\prescript{\textrm{p}}{rg}{}\Omega_s^\textrm{c}$ (of either quadrilateral or triangular shape) in the following manner:
\begin{equation}\label{eq: child_meshing}
    \bigcup_{s=1}^{\prescript{\textrm{p}}{rg}{}N^\textrm{c}} \prescript{\textrm{p}}{rg}{}\Omega_{s}^\textrm{c} = \mathcal{H}_{\mathsmaller{rg}}, \quad \prescript{\textrm{p}}{rg}{}\Omega_{i}^\textrm{c} \cap \prescript{\textrm{p}}{rg}{}\Omega_{j}^\textrm{c} = \varnothing, ~\forall i \neq j
\end{equation}
Using the spatial mesh presented above, the nonlocal stress $\bm{\sigma}$ can be expressed as:
\begin{equation}
    \bm{\sigma}(\bm{x}(\bm{\xi}_{rg}^\textrm{p}))
    = \bm{C}\left(\bm{x}(\bm{\xi}_{rg}^\textrm{p})\right) : \sum_{s=1}^{\prescript{\textrm{p}}{rg}{}N^\textrm{c}} \int_{\prescript{\textrm{p}}{rg}{}\Omega_{s}^\textrm{c}} \mathcal{K}\left(\bm{x}(\bm{\xi}_{rg}^\textrm{p}),\bm{x}^\prime\right) \bm{\epsilon}(\bm{x}^\prime) \mathrm{d}\bm{x}^\prime
\end{equation}
Particular attention is required for mesh generation when considering weakly-singular nonlocal kernel functions and Gauss-Jacobi quadrature. As discussed in \cite{pang2013gauss}, integrals with weak singularity can be numerically evaluated by Gauss-Jacobi quadrature only when the singular points coincide with the boundary; see remarks in \S\ref{sec: closing_remarks}. In this regard, the singular point (for example, $\bm{x}^\prime=\bm{x}$ in fractional-order nonlocal elasticity, see {\S}S1 of the SM) must be employed as a node in the child mesh such that it lies on the boundary of the child elements $\prescript{\textrm{p}}{rg}{}\Omega_{s}^\textrm{c}$.

\subsubsection{Numerical integration}\label{sec: 2.4.2}
In this section, we develop the approach to numerically evaluate the nonlocal stress at the various quadrature points of the parent mesh. Similar to the approach developed in \S\ref{sec: 2.3.2}, we first transform the child elements $\prescript{\textrm{p}}{rg}{}\Omega_{s}^\textrm{c}$ in 2D Cartesian coordinates into 2D isoparametric elements $\prescript{\textrm{p}}{rg}{}\hat{\Omega}_{s}^\textrm{c}$. By doing this, a given point $\bm{x}=(x,y)$ in $\prescript{\textrm{p}}{rg}{}\Omega_{s}^\textrm{c}$, can now be mapped into:
\begin{equation}
    \bm{x}(\bm{\xi}_{s}^\textrm{c}) = \sum_{d=1}^{d^\textrm{c}} \Phi_{d}^\textrm{c}(\bm{\xi}_{s}^\textrm{c}) \prescript{\textrm{p}}{rg}{}\bm{x}_{sd}^\textrm{c}
\end{equation}
where $\Phi_d^\textrm{c}$ and $d^\textrm{c}$ are the shape function and the number of nodes $\prescript{\textrm{p}}{rg}{}\bm{x}_{sd}^\textrm{c}$ in child element $\prescript{\textrm{p}}{rg}{}\Omega_{s}^\textrm{c}$, respectively. Using the Gauss quadrature principles, the integral over $\prescript{\textrm{p}}{rg}{}\Omega_{s}^\textrm{c}$ can be approximated in isoparametric elements $\prescript{\textrm{p}}{rg}{}\hat{\Omega}_{s}^\textrm{c}$, as:
\begin{equation}\label{eq: Gauss_integration_child_element}
    \int_{\prescript{\textrm{p}}{rg}{}\Omega_{s}^\textrm{c}} \mathcal{K}\left(\bm{x}(\bm{\xi}_{g}^\textrm{p}),\bm{x}^\prime\right) \bm{\epsilon}(\bm{x}^\prime) \mathrm{d}\bm{x}^\prime = 
    \sum_{h=1}^{g^\textrm{c}} w_{sh}^\textrm{c} \vert\prescript{\textrm{p}}{rg}{}\bm{J}_{s}^\textrm{c}\vert \mathcal{K}\left(\bm{x}(\bm{\xi}_{g}^\textrm{p}),\bm{x}^\prime(\bm{\xi}_{sh}^\textrm{c})\right) \bm{\epsilon}(\bm{\xi}_{sh}^\textrm{c})
\end{equation}
where $g^\textrm{c}$ is the number of Gauss points $\bm{\xi}_{sh}^\textrm{c}$ and weight $w_{sh}^\textrm{c}$ in child elements; $\vert\prescript{\textrm{p}}{rg}{}\bm{J}_{s}^\textrm{c}\vert$ is the determinant of Jacobian matrix $\prescript{\textrm{p}}{rg}{}\bm{J}_{s}^\textrm{c}$ defined by:
\begin{equation}
    \prescript{\textrm{p}}{rg}{}\bm{J}_{s}^\textrm{c} 
    = \prescript{\textrm{p}}{rg}{}\frac{\partial \bm{x}}{\partial \bm{\xi}_{s}^\textrm{c}}
    = \sum_{d=1}^{d^\textrm{c}} \frac{\partial \Phi_d^c(\bm{\xi}_{s}^\textrm{c})}{\partial \bm{\xi}_{s}^\textrm{c}} \prescript{\textrm{p}}{rg}{}\bm{x}_{sd}^\textrm{c}
    = \sum_{d=1}^{d^\textrm{c}} \prescript{\textrm{p}}{rg}{}\bm{x}_{sd}^\textrm{c} \otimes \Phi_{d,\bm{\xi}_{s}^\textrm{c}}^\textrm{p}(\bm{\xi}_{s}^\textrm{c})
\end{equation}

\paragraph*{Scale-Bridging: link to the parent domain.}
Critical insights appear from a closer analysis of the child domain-based nonlocal integral evaluation derived above. As evident from Eq.~(\ref{eq: Gauss_integration_child_element}), the determination of the stress requires an approximation of the \textit{local strain} $\bm{\epsilon}$ at the child Gauss points $\bm{\xi}_{sh}^\textrm{c}$. We first discuss the physical relevance of this numerical expression and then proceed to approximate the local strain. Note that, in the strain-driven integral elasticity approach, the (infinitesimal) strain is assumed to be the symmetric gradient of the displacement field, following local continuum assumptions. The algorithmic translation of this definition, based on the fine and coarse structural domains, indicates that the strain must be evaluated at the parent domain. In fact, it is this connection that imparts an implicit multiscale character to the overall model. Now, as the child mesh is decoupled from the parent mesh, an added numerical mechanism is required for the multiscale relay of response (here, strain).

\begin{figure}[hb!]
    \centering
    \includegraphics[width=1\linewidth]{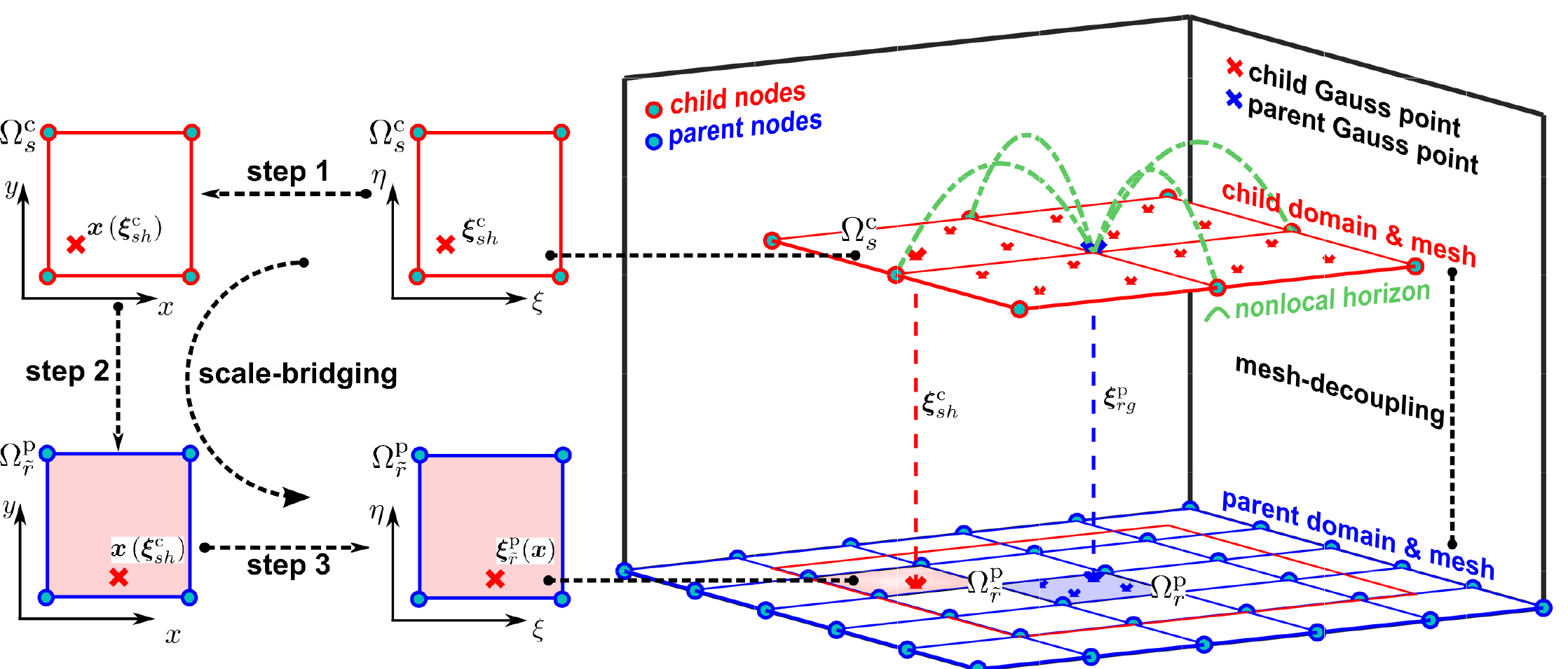}
    \caption{Schematic showing the finite element configuration including both mesh-decoupling and scale-bridging algorithms. The overall discretization of the square-shaped parent domain (or the material geometry) and the decoupled child domain (or the nonlocal horizon). The parent domain is discretized uniformly into parent meshes consisting of rectangular parent elements and nodes (see the solid blue edges and blue nodes). The rectangular-shaped nonlocal horizon for parent Gauss point $\bm{\xi}_{rg}^\textrm{p}$ (represented by the bold blue cross sign) in a given parent element $\Omega_r^\textrm{p}$ (represented by the blue patch), is discretized into individual rectangular child meshes with rectangular child elements and nodes (see the solid red edges and red nodes). The child and parent meshes are fully decoupled. The dimension of the nonlocal horizon in the four directions is outlined by the green dashed curves. The child Gauss point $\bm{\xi}_{sh}^\textrm{c}$ in the child element $\Omega_s^\textrm{c}$, highlighted by the bold red cross sign, is located at the parent element $\Omega_{\Tilde{r}}^\textrm{p}$ (represented by the red patch). The three-step scale-bridging algorithm that maps the child Gauss point $\bm{\xi}_{sh}^\textrm{c}$ to the parent point $\bm{\xi}_{\tilde{r}}^\textrm{p}$ is illustrated in the left, where child and parent element are labeled as $\Omega_s^\textrm{c}$ and $\Omega_{\tilde{r}}^\textrm{p}$, respectively. Similarly, the conversion between Euclidean and isoparametric coordinates are tracked by $x$-$y$ and $\xi$-$\eta$ in each element, respectively.}
    \label{fig: FE_configuration}
\end{figure}

We formulate a scale-bridging algorithm consisting of three different transformations to numerically implement the multiscale relay of information, described above. Without any loss of generality, we consider the approximation of $\bm{\epsilon}$ at $\bm{\xi}_{sh}^\textrm{c}$ in the parent domain level. From a high level perspective, for the scale-bridging of $\bm{\xi}_{sh}^\textrm{c}$, we first find out the parent element $\Omega_{\Tilde{r}}^\textrm{p}$ ($\neq \Omega_{{r}}^\textrm{p}$, in general; see details in step 2, below) that contains $\bm{x}(\bm{\xi}_{sh}^\textrm{c})$, and then approximate $\bm{\epsilon}(\bm{x}(\bm{\xi}_{sh}^\textrm{c}))$ within $\Omega_{\Tilde{r}}^\textrm{p}$. This scale-bridging requires the following sequential transformations:
\begin{itemize}[leftmargin=*]
    \item \textit{Step 1}: to perform the inverse isoparametric transformation from child Gauss point with isoparametric coordinates $\bm{\xi}_{sh}^\textrm{c}$ to the Euclidean coordinates $\bm{x}(\bm{\xi}_{sh}^\textrm{c})$. Note that the inverse isoparametric transformation can be achieved easily since we use linear shape functions $\Phi_d^\textrm{c}$ for geometric interpolation.
    
    \item \textit{Step 2}: to determine the specific parent element $\Omega_{\Tilde{r}}^\textrm{p}$ that contains the child Gauss point $\bm{x}(\bm{\xi}_{sh}^\textrm{c})$ (using its Cartesian coordinates). This requires a search over all parent elements because the parent element $\Omega_{\Tilde{r}}^\textrm{p}$ containing $\bm{\xi}_{sh}^\textrm{c}$, is more likely different from the reference parent element $\Omega_{{r}}^\textrm{p}$ (that is used in this \S\ref{sec: 2.4} to derive the child domain approximation; see introduction to \S\ref{sec: 2.4}). In fact, $\{\Omega_{\Tilde{r}}^\textrm{p}, \Omega_{{r}}^\textrm{p}\}$ can be identified as a pair of interacting parent elements (or equivalently, the discrete realization of the interacting pair $\{\bm{x},\bm{x}^\prime\}$ in Eq.~(\ref{eq: constitutive_relation}).
    
    \item \textit{Step 3}: to perform the forward isoparametric transformation (in the parent domain) to obtain the isoparametric coordinates of the child Gauss point $\bm{\xi}_{\Tilde{r}}^\textrm{p}$ that belongs to the parent element $\Tilde{r}$.
\end{itemize}
The three-step scale-bridging, along with underlying mesh-decoupling, are schematically illustrated in Fig.~(\ref{fig: FE_configuration}). 

Following the above-outlined strategy, $\bm{\epsilon}(\bm{x}(\bm{\xi}_{sh}^\textrm{c}))$ is approximated as:
\begin{equation}
    \bm{\epsilon}(\bm{x}(\bm{\xi}_{sh}^\textrm{c})) = \frac{1}{2}\left[ \bm{\nabla}_{\bm{x}}\bm{u}(\bm{x}(\bm{\xi}_{\Tilde{r}}^\textrm{p})) + \bm{\nabla}_{\bm{x}}^\textrm{T}\bm{u}(\bm{x}(\bm{\xi}_{\Tilde{r}}^\textrm{p})) \right]
\end{equation}
Now, analogous to Eq.~(\ref{eq: grad_delta_v}), the gradient of the displacement field at the child Gauss point $\bm{\xi}_{sh}^\textrm{c}$ can be approximated within the connected parent element $\bm{\xi}_{\Tilde{r}}^\textrm{p}$, in the following manner:
\begin{equation}
    \bm{\nabla}_{\bm{x}}\bm{u}(\bm{x}(\bm{\xi}_{\Tilde{r}}^\textrm{p}))
    = \sum_{u=1}^{d^\textrm{p}}
    \bm{u}_{\Tilde{r}u} \otimes \bm{\nabla}_{\bm{\xi}_{\Tilde{r}}^\textrm{p}} \Psi_{u}^\textrm{p}(\bm{\xi}_{\Tilde{r}}^\textrm{p}) \frac{\partial \bm{\xi}_{\Tilde{r}}^\textrm{p}}{\partial \bm{x}} 
    = \sum_{u=1}^{d^\textrm{p}} \bm{u}_{\Tilde{r}u} \otimes (\bm{J}_{\Tilde{r}}^\textrm{p})^{-T} \bm{\nabla}_{\bm{\xi}_{\Tilde{r}}^\textrm{p}} \Psi_{u}^\textrm{p}(\bm{\xi}_{\Tilde{r}}^\textrm{p}) 
\end{equation}
where $\bm{u}_{\Tilde{r}u}$ denotes DOF at nodal point $u \in \{1,2,...,d^\mathrm{p}\}$ in the parent element $\Omega_{\Tilde{r}}^\mathrm{p}$. Note that, it is the inverse isoparametric transformation that allows the mapping of the child Gauss point within the child isoparametric coordinates and the Cartesian coordinates, such that the displacement gradient at any child Gauss point can be approximated by using DOF defined in parent elements.

\subsection{Closing remarks on mesh-decoupling}
\label{sec: closing_remarks}
The mesh-decoupling and scale-bridging algorithms deserve some additional remarks. We emphasize that the mesh-decoupling essentially implies that the child mesh (and the selection of its topology) is completely independent from that of the parent mesh. Note that a physically-consistent multiscale transfer of information is achieved via both the (three-step) scale-bridging approach, and the theoretical thermodynamic consistency. The independence of the child mesh topology allows the use of highly modular unstructured child domain discretizations to account for different nonlocal kernels and irregular horizons (as it will be shown in \S\ref{sec: 5}). The weakly-singular kernels and the associated Gauss-Jacobi quadrature principle, for example, can be used only by leveraging independent child mesh since the Gauss-Jacobi principle requires the singularity point (that is the parent Gauss point $\bm{\xi}_{rg}^\textrm{p}$) to be exactly at element edges, while in comparison, if using parent mesh only (like \cite{pisano2009nonlocal}), the singular parent Gauss point $\bm{\xi}_{rg}^\textrm{p}$ will never reach the edges in parent meshes. We refer interested readers to {\S}S1 of the supplementary material (SM) for more discussion on the advantages of using the Gauss-Jacobi quadrature and the independent child mesh.

Finally, a summary of the most important steps of $\text{M}^2$-FEM in Algorithm~\ref{algm: numerical_implementation} (with a specific focus on the mesh-decoupling technique), reveals the scalable nature of the proposed concept. Evidently, the mesh-decoupling technique can be extended in either directions (coarser or finer) from the parent and/or child domains identified in Algorithm~\ref{algm: numerical_implementation}. The extension would require an additional nested set of decoupled-meshes (analogous to that used in Algorithm~\ref{algm: numerical_implementation}) and another scale-bridging to relay the relevant response information. The resulting hierarchically stacked computational model can be used to simulate highly multiscale problems. Although the need for such hierarchically stacked multiscale theories has been highlighted (see \cite{mcdowell2010perspective,chopard2014framework,ding2022multiscale}), we anticipate that the significant or prohibitive computational resources that would be required to simulate these theories by existing finite element methodologies, has prevented their development. Results presented in \S\ref{sec: 5} indicate that the proposed $\text{M}^2$-FEM could potentially enable the development and implementation of hierarchically stacked multiscale theories.

\tcbset{width=(\linewidth-4.5cm),colframe=gray!100,colback=white,height=3.60cm}
\begin{algorithm}[ht!]
\setstretch{1.05}
\DontPrintSemicolon
\SetNoFillComment

\tikzc{1pt} Approximate weak-form by evaluating all variables defined at the coarse (parent) scale\;

\tikzc{1pt} Approximate variables that are defined at fine (child) scale in weak-form\;
{\begin{tcolorbox}[size=fbox,title=Mesh-decoupling for a given parent- and child-scale combination]
\For{every element $\Omega_r^\mathrm{p}$ and corresponding quadrature points $\bm{\xi}^\mathrm{p}_r$ in the parent mesh}
{$\star$ Identify the child domain at the quadrature points $\bm{\xi}^\mathrm{p}_r$ in $\Omega^\mathrm{p}_r$\;

$\star$ Generate an independent child mesh decoupled from parent mesh\;

$\star$ Numerically evaluate variables on child mesh using discretized expressions\;

$\star$ Employ scale-bridging algorithm to relay information to parent domain\;

$\star$ Assemble all child domain approximations within weak-form integral}
\end{tcolorbox}}

\If{variables defined at another (further) finer scale exist in the weak-form}
{$\quad-$ Update the child and the finer scales as the new parent and child scales, respectively\;

$\quad-$ Repeat the parent-child mesh-decoupling approximation algorithm above, recursively}
\caption{Weak-form approximation of nested integrals in hierarchical material scales via $\text{M}^2$-FEM \label{algm: numerical_implementation}}
\end{algorithm}

\clearpage
\section{Implementation: matrix representation and assembled algebraic equations}
\label{sec: 3}
The above two sections complete the development of the finite element discretization and approximation in both the parent and child domains. Here below, we derive the matrix representation of the different finite element approximations and assemble them coherently to formulate the algebraic equations of motion. It is common knowledge that these steps are essential (and ideal) for an efficient computer implementation of the solver.

\subsection{Background: matrix representations of primary variables and expressions}
\label{sec: 3.1}
We consider the $r^{\textrm{th}}$ parent element $\Omega_r^\textrm{p}$, the displacement DOF $\bm{u}$, and the weighting function $\delta\bm{v}$ defined on the element nodes. The finite element approximation of these quantities can be expressed as the following vectors:
\begin{subequations}
\begin{align}
    \{\bm{u}\}_r^\textrm{p} &=
    \{\{\bm{u}_x\}_r^\mathrm{p}, \{\bm{u}_y\}_r^\mathrm{p}\}^\mathrm{T} = 
    \{\{\bm{u}_{x1},\bm{u}_{x2},...,\bm{u}_{xd^\mathrm{p}}\}, \{\bm{u}_{y1},\bm{u}_{y2},...,\bm{u}_{yd^\mathrm{p}}\}\}^\mathrm{T} \\
    \{\delta\bm{v}\}_r^\mathrm{p} &= 
    \{\{\delta\bm{v}_x\}_r^\mathrm{p}, \{\delta\bm{v}_y\}_r^\mathrm{p}\}^\textrm{T} =
    \{\{\delta\bm{v}_{x1},\delta\bm{v}_{x2},...,\delta\bm{v}_{xd^\textrm{p}}\},
    \{\delta\bm{v}_{y1},\delta\bm{v}_{y2},...,\delta\bm{v}_{yd^\textrm{p}}\}\}^\textrm{T}
\end{align}
\end{subequations}
in which $\bm{u}_{xr}$ and $\bm{u}_{yr}$ (or, $\delta\bm{v}_{xr}$ and $\delta\bm{v}_{yr}$) are the previously defined solution vectors in size $k \times 1$, that represent the nodal DOF (or, weighting functions) in the $x$ and $y$ direction, respectively. In that regard, the DOF $\{\bm{u}_x\}_r^\textrm{p}$ and $\{\bm{u}_y\}_r^\textrm{p}$ (or, $\{\delta\bm{v}_x\}_r^\textrm{p}$ and $\{\delta\bm{v}_y\}_r^\textrm{p}$) in a given parent element are all row vectors of size $m^\textrm{p} \times 1$. Consequently, the length of vectors $\{\bm{u}\}_r^\textrm{p}$ and $\{\delta\bm{v}\}_r^\textrm{p}$ is $2m^\textrm{p}$. Similarly, the DOF and weighting functions on the boundary element $\Gamma_r^\textrm{t}$, both in size $2m^\textrm{t} \times 1$, can be expressed as:
\begin{subequations}
\begin{align}
    \{\bm{u}\}_r^\textrm{t} &=
    \{\{\bm{u}_x\}_r^\textrm{t}, \{\bm{u}_y\}_r^\textrm{t}\}^\textrm{T} = 
    \{\{\bm{u}_{x1},\bm{u}_{x2},...,\bm{u}_{xd^\textrm{t}}\}, \{\bm{u}_{y1},\bm{u}_{y2},...,\bm{u}_{yd^\textrm{t}}\}\}^\textrm{T} \\
    \{\delta\bm{v}\}_r^\textrm{t} &= 
    \{\{\delta\bm{v}_x\}_r^\textrm{t}, \{\delta\bm{v}_y\}_r^\textrm{t}\}^\textrm{T} =
    \{\{\delta\bm{v}_{x1},\delta\bm{v}_{x2},...,\delta\bm{v}_{xd^\textrm{t}}\},
    \{\delta\bm{v}_{y1},\delta\bm{v}_{y2},...,\delta\bm{v}_{yd^\textrm{t}}\}\}^\textrm{T}
\end{align}
\end{subequations}

We retrace the above steps for developing matrix representations of the shape functions $\bm{\Phi}^\textrm{p}$ and $\bm{\Psi}^\textrm{p}$ in parent elements, that perform the isogeometric transformation and the interpolation of solution variables (also weighting functions, see Eq.~(\ref{eq: u_and_delta_v})), respectively. Following their definitions, we have:
\begin{subequations}
\begin{align}
    \{\bm{\Phi}^\textrm{p}(\bm{\xi}_r^\textrm{p})\} &= \{\Phi_1(\bm{\xi}_r^\textrm{p}),\Phi_2(\bm{\xi}_r^\textrm{p}),...,\Phi_{d^\textrm{p}}(\bm{\xi}_r^\textrm{p})\}^\textrm{T} \\
    \{\bm{\Psi}^\textrm{p}(\bm{\xi}_r^\textrm{p})\} &= \{\Psi_1^\textrm{p}(\bm{\xi}_r^\textrm{p}),\Psi_2^\textrm{p}(\bm{\xi}_r^\textrm{p}),...,\Psi_{m^\textrm{p}}^\textrm{p}(\bm{\xi}_r^\textrm{p})\}^\textrm{T}
\end{align}
\end{subequations}
As we have clarified in \S\ref{sec: 2.2}, depending on the interpolation, the two sets of shape functions can differ in both the definitions and dimensions (note that $\{\bm{\Phi}\}^\textrm{p}$ is in size $d^\textrm{p} \times 1$ and $\{\bm{\Psi}^\textrm{p}\}$ is in size $m^\textrm{p} \times 1$). In analogy with the above formulation, shape functions $\bm{\Phi}^\textrm{t}$ and $\bm{\Psi}^\textrm{t}$ defined on boundary elements $\Gamma_r^\textrm{t}$ can be also expressed as:
\begin{subequations}
\begin{align}
    \{\bm{\Phi}^\textrm{t}(\bm{\xi}_r^\textrm{t})\} &= \{\Phi_1(\bm{\xi}_r^\textrm{t}),\Phi_2(\bm{\xi}_r^\textrm{t}),...,\Phi_{d^\textrm{t}}(\bm{\xi}_r^\textrm{t})\}^\textrm{T} \\
    \{\bm{\Psi}^\textrm{t}(\bm{\xi}_r^\textrm{t})\} &= \{\Psi_1^\textrm{t}(\bm{\xi}_r^\textrm{t}),\Psi_2^\textrm{t}(\bm{\xi}_r^\textrm{t}),...,\Psi_{m^\textrm{t}}^\textrm{t}(\bm{\xi}_r^\textrm{t})\}^\textrm{T}
\end{align}
\end{subequations}

We recast the original tensorial weak-form governing equations into the Voigt matrix form. For nonlocal plane-strain problems with isotropic and homogeneous elastic media, the second-order strain tensor, fourth-order stiffness tensor, and the second-order nonlocal stress tensor can be recast into matrices as:
\begin{equation}
    [\bm{C}]_{3 \times 3} =
    \begin{bmatrix}
        2\mu+\lambda & \lambda & 0 \\
        \lambda & 2\mu+\lambda & 0 \\
        0 & 0 & 2\mu \\
    \end{bmatrix}, \quad
    \{\bm{\sigma}\}_{3 \times 1} =
    \begin{bmatrix}
        \sigma_{xx} \\
        \sigma_{yy} \\
        \sigma_{xy} 
    \end{bmatrix}, \quad
    \{\bm{\epsilon}\}_{3 \times 1} =
    \begin{bmatrix}
        \epsilon_{xx} \\
        \epsilon_{yy} \\
        \epsilon_{xy} 
    \end{bmatrix}
\end{equation}
such that the nonlocal constitutive relation in Eq.~(\ref{eq: constitutive_relation}) can be expressed as:
\begin{equation}\label{eq: constitutive_relation_matrix}
    \{\bm{\sigma}\} = [\bm{C}]\int_{\mathcal{H}}\mathcal{K}(\bm{x},\bm{x}^\prime)\{\bm{\epsilon}(\bm{x}^\prime)\} \mathrm{d}\bm{x}^\prime
\end{equation}
By using this technique, the double dot product $\bm{\sigma} : \bm{\nabla}_{\bm{x}}\delta\bm{v}$ in the weak-form formulation (see Eq.~(\ref{eq: weak_governing_eqn})) can be also recast into the following matrix representation:
\begin{equation}
    \bm{\sigma} : \bm{\nabla}_{\bm{x}}\delta\bm{v} = \{\bm{\sigma}\}^\textrm{T}\{\bm{\nabla}_{\bm{x}}^{\textrm{sym}}\delta\bm{v}\}_{3 \times 1} = 
    \{\bm{\sigma}\}^\textrm{T}
    \begin{bmatrix}
        \delta\bm{v}_{x,x} &   \delta\bm{v}_{y,y} &
        \dfrac{1}{2}\left(\delta\bm{v}_{x,y} + \delta\bm{v}_{y,x} \right)
    \end{bmatrix}^\textrm{T}
\end{equation}
where $\bm{\nabla}_{\bm{x}}^{\textrm{sym}}$ denotes the symmetric gradient in Voigt notation. By using the above framework, we recast the $\text{M}^2$-FEM algorithm into an algebraic representation that can be easily translated into a scientific computing program.

\subsection{Matrix representation of numerically evaluated approximations}
\label{sec: 3.2}
We first evaluate the terms containing the weighting functions $\delta\bm{v}$. Specifically, we start with the finite element approximation of the weighting functions $\delta\bm{v}(\bm{\xi}_{rg}^\textrm{p})$ at Gauss points $\bm{\xi}_{rg}^\textrm{p}$. By combining the numerical scheme in Eq.~(\ref{eq: u_and_delta_v}) and the vectorial form of shape functions $\bm{\Psi}^\textrm{p}$, we obtain the matrix representation of $\delta\bm{v}(\bm{\xi}_{rg}^\textrm{p})$ as follows:
\begin{equation}\label{eq: delta_v_vector}
    \{\delta\bm{v}(\bm{\xi}_{rg}^\textrm{p})\} = 
    \{\bm{L}\}_{1 \times 2m^\textrm{p}}^{(\bm{\Psi}^\textrm{p})}\{\delta\bm{v}\}_r =
    \{\bm{\Psi}^\textrm{p}(\bm{\xi}_r^\textrm{p}),\bm{\Psi}^\textrm{p}(\bm{\xi}_r^\textrm{p})\}^\textrm{T}\{\delta\bm{v}\}_r
\end{equation}
Note that we use $\{\bm{\Psi}^\textrm{p}(\bm{\xi}_r^\textrm{p})\}$ twice in the coefficient matrix $\{\bm{L}\}^{(\bm{\Psi}^\textrm{p})}$ to account for the approximation of both ${\delta\bm{v}_x}_r$ and ${\delta\bm{v}_y}_r$. The other weighting function related term in the parent element, which is the gradient of weighting function $\bm{\nabla}_{\bm{x}} \delta\bm{v}(\bm{\xi}_{rg}^\textrm{p})$ in Eq.~(\ref{eq: grad_delta_v}), can be replaced by its vectorial form $\{\bm{\nabla}_{\bm{x}}^{\textrm{sym}}\delta\bm{v}(\bm{\xi}_{rg}^\textrm{p})\}$:
\begin{equation}\label{eq: sym_grad_delta_v_vector}
    \{\bm{\nabla}_{\bm{x}}^{\textrm{sym}}\delta\bm{v}(\bm{\xi}_{rg}^\textrm{p})\} =
    [\bm{G}]_{3 \times 2m^\textrm{p}}^{(\delta\bm{v})}\{\delta\bm{v}\}_r =
    \frac{1}{2}
    \begin{bmatrix}
        \{\bm{G}_{rx}^\textrm{p}(\bm{\xi}_{rg}^\textrm{p})\}^\textrm{T} & \{\bm{0}\}_{1 \times m^\textrm{p}} \\[4pt]
        \{\bm{0}\}_{1 \times m^\textrm{p}} & \{\bm{G}_{ry}^\textrm{p}(\bm{\xi}_{rg}^\textrm{p})\}^\textrm{T} \\[4pt]
        \{\bm{G}_{rx}^\textrm{p}(\bm{\xi}_{rg}^\textrm{p})\}^\textrm{T} &
        \{\bm{G}_{ry}^\textrm{p}(\bm{\xi}_{rg}^\textrm{p})\}^\textrm{T} 
    \end{bmatrix}
    \{\delta\bm{v}\}_r
\end{equation}
where $[\bm{G}]_{3 \times 2m^\textrm{p}}^{(\delta\bm{v})}$ is the coefficient matrix in size $3 \times 2m^\textrm{p}$; $\{\bm{0}\}_{1 \times m^\textrm{p}}$ is the zero vector in size $1 \times m^\textrm{p}$; $\{\bm{G}_{rx}^\textrm{p}(\bm{\xi}_{rg}^\textrm{p})\}$ and $\{\bm{G}_{rx}^\textrm{p}(\bm{\xi}_{rg}^\textrm{p})\}$ are two coefficient vectors in size $m^\textrm{p} \times 1$ that represent the first-order spatial gradient along $x$ and $y$ directions. The two gradient coefficient vectors can be evaluated by using the matrix $[\bm{G}_{r}^\textrm{p}(\bm{\xi}_{rg}^\textrm{p})]$ (see Eq.~(\ref{eq: grad_delta_v})):
\begin{equation}\label{eq: grad_matrix}
    [\{\bm{G}_{rx}^\textrm{p}(\bm{\xi}_{rg}^\textrm{p})\}, \{\bm{G}_{ry}^\textrm{p}(\bm{\xi}_{rg}^\textrm{p})\}] =
    [\bm{G}_{r}^\textrm{p}(\bm{\xi}_{rg}^\textrm{p})]_{m^\textrm{p} \times 2} = 
    [\bm{\nabla}_{\bm{\xi}_{r}^\textrm{p}} 
        \Psi_{v}^\textrm{p}(\bm{\xi}_{rg}^\textrm{p}) (\bm{J}_{rg}^\textrm{p})^{-1}] = 
    \begin{bmatrix}
        \bm{\nabla}_{\bm{\xi}_{r}^\textrm{p}} 
        \Psi_{1}^\textrm{p}(\bm{\xi}_{rg}^\textrm{p}) (\bm{J}_{rg}^\textrm{p})^{-1} \\[4pt]
        \bm{\nabla}_{\bm{\xi}_{r}^\textrm{p}} 
        \Psi_{2}^\textrm{p}(\bm{\xi}_{rg}^\textrm{p}) (\bm{J}_{rg}^\textrm{p})^{-1} \\[4pt]
        \vdots \\[4pt]
        \bm{\nabla}_{\bm{\xi}_{r}^\textrm{p}} 
        \Psi_{m^{\textrm{p}}}^\textrm{p}(\bm{\xi}_{rg}^\textrm{p}) (\bm{J}_{rg}^\textrm{p})^{-1}
    \end{bmatrix}
\end{equation}
Finally, the matrix form of the approximated $\delta\bm{v}$ on boundary elements $\Gamma_{r}^\textrm{t}$ can be readily given by:
\begin{equation}
    \{\delta\bm{v}(\bm{\xi}_{rg}^\textrm{t})\} = 
    \{\bm{L}(\bm{\xi}_{rg}^\textrm{t})\}_{1 \times 2m^\textrm{t}}^{(\bm{\Psi}^\textrm{t})} \{\delta\bm{v}\}_r^\textrm{t} =
    \{\bm{\Psi}^\textrm{t}(\bm{\xi}_{rg}^\textrm{t}),\bm{\Psi}^\textrm{t}(\bm{\xi}_{rg}^\textrm{t})\}^\textrm{T}\{\delta\bm{v}\}_r^\textrm{t}
\end{equation}

We then derive the matrix representation for terms containing $\bm{u}$. Similar to Eq.~(\ref{eq: delta_v_vector}), the matrix representation of $\Ddot{\bm{u}}(\bm{\xi}_{rg}^\textrm{p})$ shown in Eq.~(\ref{eq: u_and_delta_v}) can be obtained immediately as:
\begin{equation}
    \{\Ddot{\bm{u}}(\bm{\xi}_{rg}^\textrm{p})\} = 
    \{\bm{L}(\bm{\xi}_{rg}^\textrm{p})\}_{1 \times 2m^\textrm{p}}^{(\bm{\Psi}^\textrm{p})}\{\Ddot{\bm{u}}\}_r
\end{equation}
The matrix formulation with respect to the other solution variable related term $\bm{\sigma}$ is more complicated due to the nonlocal integration. To formulate the nonlocal stress $\bm{\sigma}$, we start with the local strain tensor $\bm{\epsilon}(\bm{x}(\bm{\xi}_{sh}^\textrm{c}))$ defined at the child Gauss point $\bm{\xi}_{sh}^\textrm{c}$. In analogy with $\{\bm{\nabla}_{\bm{x}}^{\textrm{sym}}\delta\bm{v}\}$ in Eq.~(\ref{eq: sym_grad_delta_v_vector}), the second-order strain tensor $\bm{\epsilon}$ can be also reformulated accordingly into the following Voigt form:
\begin{equation}\label{eq: strain_vector}
     \{\bm{\epsilon}(\bm{x}(\bm{\xi}_{sh}^\textrm{c}))\} =
     [\bm{G}(\bm{\xi}_{sh}^\textrm{c})]_{3 \times 2m^\textrm{p}}^{(\bm{u})}\{\bm{u}\}_r =
     \frac{1}{2}
    \begin{bmatrix}
         \{\bm{G}_{\Tilde{r}x}^\textrm{p}(\bm{\xi}_{\Tilde{r}g}^\textrm{p})\}^\textrm{T} & \{\bm{0}\}_{1 \times m^\textrm{p}} \\[4pt]
         \{\bm{0}\}_{1 \times m^\textrm{p}} & \{\bm{G}_{\Tilde{r}y}^\textrm{p}(\bm{\xi}_{\Tilde{r}g}^\textrm{p})\}^\textrm{T} \\[4pt]
         \{\bm{G}_{\Tilde{r}y}^\textrm{p}(\bm{\xi}_{\Tilde{r}g}^\textrm{p})\}^\textrm{T} &
         \{\bm{G}_{\Tilde{r}x}^\textrm{p}(\bm{\xi}_{\Tilde{r}g}^\textrm{p})\}^\textrm{T}
    \end{bmatrix}
    \{\bm{u}\}_{\Tilde{r}}
\end{equation}
in which the $\{\bm{G}_{\Tilde{r}x}^\textrm{p}(\bm{\xi}_{\Tilde{r}g}^\textrm{p})\}$ and $\{\bm{G}_{\Tilde{r}y}^\textrm{p}(\bm{\xi}_{\Tilde{r}g}^\textrm{p})\}$ share the same definition with their counterparts, namely the $\{\bm{G}_{rx}^\textrm{p}(\bm{\xi}_{rg}^\textrm{p})\}$ and $\{\bm{G}_{rx}^\textrm{p}(\bm{\xi}_{rg}^\textrm{p})\}$, that are defined in Eq.~(\ref{eq: grad_matrix}). Note that $\bm{\xi}_{sh}^\textrm{c}$ and $\bm{\xi}_{\Tilde{r}g}^\textrm{p}$ are essentially the same child Gauss point, but represented in the child and parent mesh, respectively (see the discussion in \S\ref{sec: 2.4.2}). Further substituting Eq.~(\ref{eq: strain_vector}) into the nonlocal integral, we obtain the matrix representation of the nonlocal stress $\bm{\sigma}$:
\begin{equation}\label{eq: stress_vector}
    \bm{\sigma}(\bm{x}(\bm{\xi}_{rg}^\textrm{p})) 
    = [\bm{E}(\bm{\xi}_{rg}^\textrm{p})]_{3 \times 2kN^\textrm{s}}\{\bm{u}\}_{2kN^\textrm{s} \times 1}
\end{equation}
with the coefficient matrix $[\bm{E}(\bm{\xi}_{rg}^\textrm{p})]$ in size $3 \times 2kN^\textrm{s}$ is calculated by:
\begin{equation}
    [\bm{E}(\bm{\xi}_{rg}^\textrm{p})] = 
    [\bm{C}(\bm{x}(\bm{\xi}_{rg}^\textrm{p}))] \left[\sum_{s=1}^{\prescript{\textrm{p}}{rg}{}N^\textrm{c}}\sum_{h=1}^{n^\textrm{c}} \omega_{sh}^\textrm{c}|\prescript{\textrm{p}}{rg}{}\bm{J}_{s}^\textrm{c}| \mathcal{K}\left(\bm{x}(\bm{\xi}_{g}^\textrm{p}),\bm{x}^\prime(\bm{\xi}_{sh}^\textrm{c})\right)     [\bm{G}(\bm{\xi}_{sh}^\textrm{c})]^{(\bm{u})}\right]
\end{equation}
Note that, due to the nonlocal integration, the coefficient matrix $[\bm{E}(\bm{\xi}_{rg}^\textrm{p})]$ for the nonlocal stress $\bm{\sigma}$ is not in shape $3 \times 2m^\textrm{p}$ any more; instead, it requires displacement information from all the nodes within the nonlocal horizon at the current parent Gauss point $\bm{\xi}_{rg}^\textrm{p}$. Depending on the size of nonlocal horizon, the total number of nodes interacting with the current parent Gauss point varies. For this purpose, in Eq.~(\ref{eq: stress_vector}) we define the coefficient matrix to be in size $3 \times 2kN^\textrm{s}$ such that it allows including potential nonlocal interactions from all the $N^\textrm{s}$ nodes in the parent domain. Accordingly, we use the global DOF vector $\{\bm{u}\}$ (which is in size $2kN^\textrm{s} \times 1$) to denote DOF at all the $N^\textrm{s}$ nodes.

\subsection{Assembly and derivation of algebraic equations of motion}
\label{sec: 3.3}
The matrix representation of the element-wise weak-form integration now can be derived readily. By assembling all matrices into Eqs.~(\ref{eq: Gauss_integration_parent_element_1},\ref{eq: Gauss_integration_parent_boun_1}), the weak-form integration can be reorganized into the following matrix forms:
\begin{subequations}
\begin{align}
    \int_{\Omega_r^\textrm{p}} \bm{\sigma} : \bm{\nabla}_{\bm{x}}\delta\bm{v} \mathrm{d}\Omega_r^\textrm{p}
    &= \{\bm{u}\}^\textrm{T}[\bm{K}_r^\textrm{p}]_{2kN^\textrm{s} \times 2m^\textrm{p}}\{\delta\bm{v}\}_r^\textrm{p} \\
    \int_{\Omega_r^\textrm{p}}\rho(\bm{f}-\Ddot{\bm{u}}) \delta\bm{v} \mathrm{d}\Omega_r^\textrm{p} 
    &= \{\bm{F}_r^\textrm{p}\}_{1 \times 2m^\textrm{p}}\{\delta\bm{v}\}_r^\textrm{p}
    - \left[\{\Ddot{\bm{u}}\}_r^\textrm{p}\right]^\textrm{T}[\bm{M}_r^\textrm{p}]_{2m^\textrm{p} \times 2m^\textrm{p}}\{\delta\bm{v}\}_r^\textrm{p} \\
    \int_{\Gamma_r^\textrm{t}} \overline{\bm{t}}\delta\bm{v} \mathrm{d}\Gamma_r^\textrm{t}  
    &= \{\bm{F}_r^\textrm{t}\}_{1 \times 2m^\textrm{t}}\{\delta\bm{v}\}_r^\textrm{t}
\end{align}
\end{subequations}
where $[\bm{K}_r^\textrm{p}]$, $[\bm{M}_r^\textrm{p}]$, $\{\bm{F}_r^\textrm{p}\}$, and $\{\bm{F}_r^\textrm{t}\}$ are coefficient matrices and vectors that are defined, respectively, as:
\begin{subequations}
\begin{align}
    [\bm{K}_r^\textrm{p}] &= \sum_{g=1}^{g^\textrm{p}} {w}_{rg}^\textrm{p} \vert\bm{J}_r^\textrm{p}(\bm{\xi}_{rg}^\textrm{p})\vert 
    [\bm{E}(\bm{\xi}_{rg}^\textrm{p})]^\textrm{T}
    [\bm{G}(\bm{\xi}_{rg}^\textrm{p})]^{(\delta\bm{v})} \\
    [\bm{M}_r^\textrm{p}] &= \sum_{g=1}^{g^\textrm{p}} w_{rg}^\textrm{p} \vert\bm{J}_r^\textrm{p}(\bm{\xi}_{rg}^\textrm{p})\vert
    \rho(\bm{\xi}_{rg}^\textrm{p})
    \left[\{\bm{L}(\bm{\xi}_{rg}^\textrm{p})\}^{(\bm{\Psi}^\textrm{p})}\right]^\textrm{T}
    \{\bm{L}(\bm{\xi}_{rg}^\textrm{p})\}^{(\bm{\Psi}^\textrm{p})} \\
    \{\bm{F}_r^\textrm{p}\} &= \sum_{g=1}^{g^\textrm{p}} w_{rg}^\textrm{p} \vert\bm{J}_r^\textrm{p}(\bm{\xi}_{rg}^\textrm{p})\vert
    \rho(\bm{\xi}_{rg}^\textrm{p}) \bm{f}(\bm{\xi}_{rg}^\textrm{p})
    \{\bm{L}(\bm{\xi}_{rg}^\textrm{p})\}^{(\bm{\Psi}^\textrm{p})} \\
    \{\bm{F}_r^\textrm{t}\} &= \sum_{g=1}^{g^\textrm{t}} w_{rg}^\textrm{t}\vert\bm{J}_r^\textrm{t}((\bm{\xi}_{rg}^\textrm{t}) )\vert \overline{\bm{t}}(\bm{\xi}_{rg}^\textrm{t}) 
    \{\bm{L}(\bm{\xi}_{rg}^\textrm{t})\}^{(\bm{\Psi}^\textrm{t})}
\end{align}
\end{subequations}
Iterating the above process for integration over all the parent and boundary elements, we obtain the final matrix representation of the finite element approximation to the weak-form governing equations in Eq.~(\ref{eq: weak_governing_eqn}):
\begin{equation}\label{eq: weak_form_eqn_matrix}
    \{\Ddot{\bm{u}}\}^\textrm{T} [\bm{M}]_{2kN^\textrm{s} \times 2kN^\textrm{s}} \{\delta\bm{v}\} +
    \{\bm{u}\}^\textrm{T} [\bm{K}]_{2kN^\textrm{s} \times 2kN^\textrm{s}} \{\delta\bm{v}\} =
    \{\bm{F}\}_{1 \times 2kN^\textrm{s}}\{\delta\bm{v}\}
\end{equation}
by assembling all the coefficient matrices and vectors $[\bm{K}_r^\textrm{p}]$, $[\bm{M}_r^\textrm{p}]$, $\{\bm{F}_r^\textrm{p}\}$, and $\{\bm{F}_r^\textrm{t}\}$ in each elements together. To guarantee that Eq.~(\ref{eq: weak_form_eqn_matrix}) holds under any variation of weighting function $\{\delta\bm{v}\}$, we must have:
\begin{equation}\label{eq: FE_linear_system}
    [\bm{M}]^\textrm{T}\{\Ddot{\bm{u}}\} + [\bm{K}]^\textrm{T}\{\bm{u}\} = \{\bm{F}\}^\textrm{T}
\end{equation}
where $[\bm{M}]^\textrm{T}$, $[\bm{K}]^\textrm{T}$, and $\{\bm{F}\}^\textrm{T}$ are the final mass matrix, the stiffness matrix, and the force vector, respectively (note that the transpose operators are applied for convenience). Equation~(\ref{eq: FE_linear_system}) is the finally approximated algebraic linear system equivalent of the integral nonlocal elasticity problem presented in Eq.~(\ref{eq: governing_eqn_parent_domain}).

\section{Validation and convergence studies}
\label{sec: 4}
In this section, we investigate the validity and convergence of the developed algorithm by simulating a set of 2D nonlocal plane strain problems. A common parent domain, with size $L_x \times L_y = L\times L = 1\textrm{m} \times 1\textrm{m}$ and Lam\'e's parameter $\mu = 1\textrm{Pa}$ and $\lambda = 1\textrm{Pa}$, is adopted for all the simulations. The material is chosen to be isotropic in nature following the assumptions of Eringen's theory \cite{eringen1984theory,batra159misuse}. In order to investigate the performance of $\text{M}^2$-FEM, we define and perform two parametric case studies by varying the definition of the kernel and the size of the horizon.

\subsection{Definitions of case studies}
\label{sec: 4.1}
Two case studies are adopted:
\begin{itemize}[leftmargin=*]
    \item \textit{Case study 1 - Non-singular bi-exponential kernel} adopted in \cite{eringen1984theory} is assumed:
    \begin{equation}\label{eq: kernel_example_1}
        \mathcal{K}(x,y,x^\prime,y^\prime,\tau) = 
            \dfrac{1}{\pi\tau} \mathrm{exp}\left(-\dfrac{|x-x^\prime|^2}{\tau}\right) \mathrm{exp}\left(-\dfrac{|y-y^\prime|^2}{\tau}\right),\quad (x^\prime,y^\prime) \in \mathcal{H}(x,y)
    \end{equation}
    where $\tau$ is the nonlocal strength parameter that controls the decay of the bi-exponential nonlocal kernel function.

    \item \textit{Case study 2 - Weakly-singular power-law kernel} adopted in fractional-order nonlocal theories \cite{ding2022multiscale,ding2022multiscale2}, is assumed:
    \begin{equation}\label{eq: kernel_example_2}
        \mathcal{K}(x,y,x^\prime,y^\prime,\alpha) = \frac{1}{\Gamma^2(1-\alpha)}\frac{1}{|x-x^\prime|^{\alpha}|y-y^\prime|^{\alpha}},\quad (x^\prime,y^\prime) \in \mathcal{H}(x,y)
    \end{equation}
    where $\alpha \in (0,1)$ controls the spatial decay of the kernel and is called as 'fractional-order', and $\Gamma(\cdot)$ denotes the Gamma function \cite{das2011functional}. Note that the power-law kernel is (weakly) singular along lines $x^\prime = x$ and $y^\prime = y$.
\end{itemize}
Figure~(S1) in the SM illustrates the difference between the above kernel functions. Note that, both the kernel functions defined above are positive-definite and symmetric such that the nonlocal model is well-posed (see \S\ref{sec: Brief_theory} or \cite{polizzotto2001nonlocal}). In each case study, the nonlocal horizon is assumed to be isotropic and symmetric in shape such that, at a given point $(x,y)$ in the parent domain, the 2D nonlocal horizon can be expressed as:
\begin{equation}
\label{eq: nonlocal_horizon_verification}
    \mathcal{H}(x,y) = [\textrm{min}(x-l^\textrm{nl},0),\textrm{max}(x+l^\textrm{nl},L)] \times [\textrm{min}(y-l^\textrm{nl},0),\textrm{max}(y+l^\textrm{nl},L)]
\end{equation}
where $l^\textrm{nl}$ denotes the length-scale in the child domain. Note that the minimum and maximum operators are required to contain the nonlocal horizon within the material boundaries.

For the present analyses, the parent and child domains are assumed to be squares for two specific reasons. First, this choice enables the derivation of closed-form analytical expressions for specific loading conditions which are utilized for validation. Second, this choice allows the definition of commonly used computational resolution measures which provide clear insights on the convergence characteristics of the integral solver. Next, the parent domain is uniformly discretized with elements of size $\Delta_x^\textrm{p} = \Delta_y^\textrm{p} = \Delta^\textrm{p}$, such that the parent mesh number is $M_x^\textrm{p} = M_y^\textrm{p} = M^\textrm{p} = L/\Delta^\textrm{p}$. Unlike the parent mesh, the child mesh is unstructured due to the truncation of nonlocal horizons at boundaries, and also the possible assumption of weakly-singular nonlocal kernels (recall from \S\ref{sec: 3.2}, in this case the Gauss quadrature points must coincide with the nodes). We have discussed the child domain discretization strategy in \S1 of the SM, in the interest of keeping the manuscript concise. Given the non-uniform child mesh, we use the spatial average of the child element size, denoted as $\Delta^\textrm{c}$, to represent the child element size. The ratio of the parent and child element sizes is denoted as $r^{\textrm{p-c}} (=\Delta^\textrm{p}/\Delta^\textrm{c})$, and called as relative resolution.

\subsection{Validation studies}
\label{sec: 4.2}
We validate the $\text{M}^2$-FEM by comparing the results of the solver against the exact solution (the analytical expressions of the displacement field) of the defined set (\S\ref{sec: 4.1}) of 2D plane strain integral elasticity problems with spatially-varying loads. Since, analytical solutions of differ-integral partial differential equations are not derived readily (even for very simple loading conditions), we adopt a different strategy, as discussed in the following. First, we assume the 2D displacement field to obey the following analytical variation:
\begin{equation}\label{eq: displacement_verification}
    \bm{u}(x,y) \equiv \{u_x(x,y),u_y(x,y)\} = \{(-x^2+L x)(-y^2+L y),0\}
\end{equation}
The assumed field satisfies homogeneous displacement boundary conditions, and the bi-parabolic $u_x$ displacement achieves its maximum value $\Bar{u}_x^{\textrm{ref}}=0.0625$[m] at $(L/2,L/2)$. The above reference displacement field is substituted within the nonlocal model in \S\ref{sec: Brief_theory} to obtain the analytical expressions of the body loads ($f_x$ and $f_y$); see SM \S2. Finally, the derived loads are used as input to $\text{M}^2$-FEM to simulate the case studies and evaluate the accuracy of the numerical results (against the reference solution in Eq.~(\ref{eq: displacement_verification})).

\begin{figure}[ht!]
    \centering
    \includegraphics[width=\linewidth]{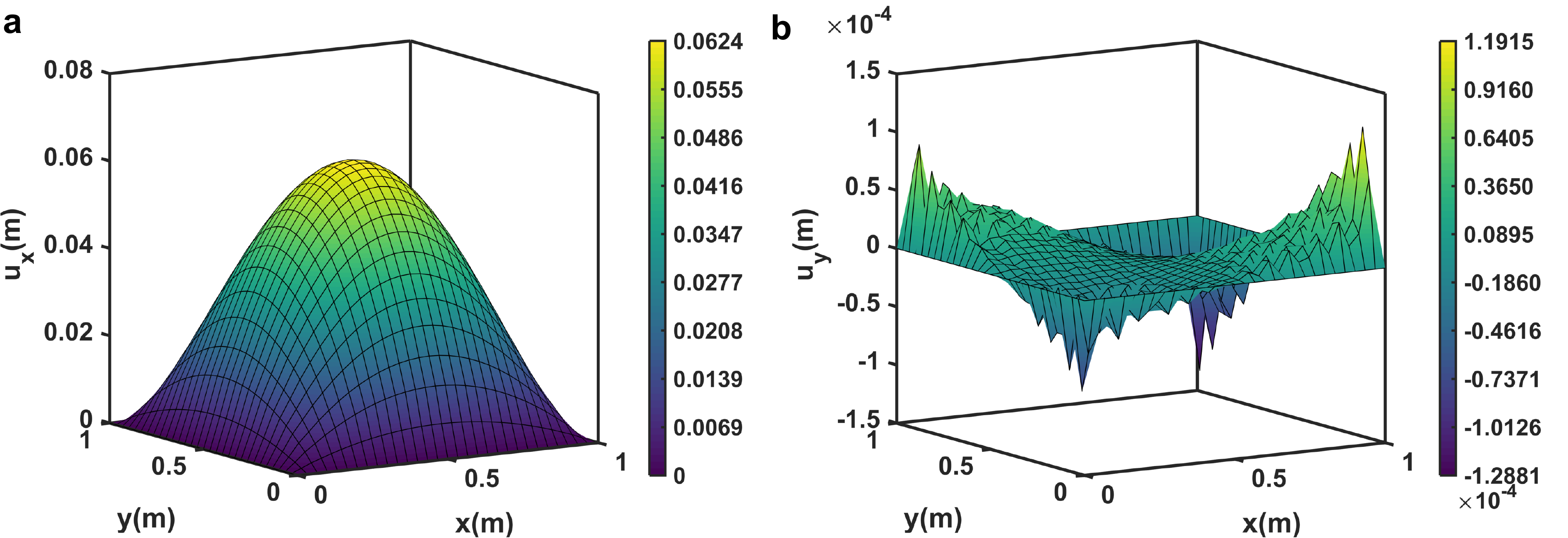}
    \caption{$\text{M}^2$-FEM simulation of the nonlocal elastic configuration with non-singular bi-exponential kernel in \textit{Case study 1}. The nonlocal response is presented in terms of the displacement field along the \textbf{a} $x-$direction ($u_y$) and \textbf{b} $x-$direction ($u_y$). The nonlocal model parameters are set to $\tau_x=\tau_y=0.002$ and $l_x=l_y=0.5$, and the mesh parameters are set to $\Delta_x^\textrm{p} = \Delta_y^\textrm{p} = 1/30$ and $\Delta_x^\textrm{c} = \Delta_y^\textrm{c} = 1/30$. Both the weak-form integral in the parent domain and the bi-exponential nonlocal integral in the child domain are evaluated numerically using Gauss-Legendre quadrature. The second-order Lagrange element is used for interpolation.}
    \label{fig: verification_case_1}
\end{figure}

\begin{figure}[htbp]
    \centering
    \includegraphics[width=\linewidth]{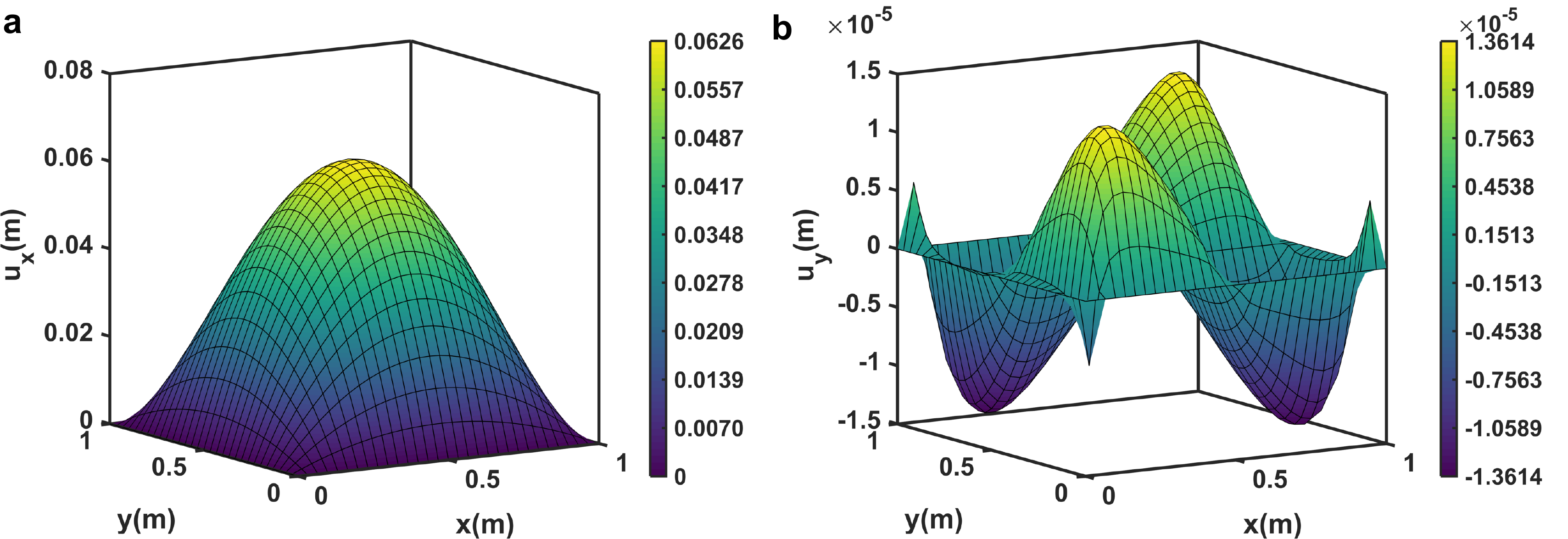}
    \caption{$\text{M}^2$-FEM simulation of the nonlocal elastic configuration with singular (fractional-order) power-law kernel in \textit{Case study 2}. The nonlocal response is presented in terms of the displacement field along the \textbf{a} $x-$direction ($u_y$) and \textbf{b} $x-$direction ($u_y$). The nonlocal model parameters are set as $\alpha_x=\alpha_y=0.5$ and $l_x=l_y=0.5$, and the mesh parameters are set as $\Delta_x^\textrm{p} = \Delta_y^\textrm{p} = 1/30$ and $\Delta_x^\textrm{c} = \Delta_y^\textrm{c} = 1/30$. The weak-form integral is numerically evaluated by the Gauss-Legendre quadrature and the weakly-singular nonlocal integral is evaluated by using Gauss-Jacobi quadrature. The first-order Lagrange element is for interpolation.}
    \label{fig: verification_case_2}
\end{figure}

Using the above strategy, the error noted in the simulation of a set of parametric configurations, obtained by varying the nonlocal horizon ($l^\textrm{nl}$) and the kernel decay-rate parameter ($\tau$ for case study  and $\alpha$ for case study 2), are evaluated. Additionally, for a comprehensive analysis, the resolution of either (a) parent mesh or (b) child mesh is varied. The results are are presented in Table~\ref{tab: Exp_accuracy} and Table~\ref{tab: PL_accuracy} (for the case study 1 and case study 2, respectively) in terms of the simulated maximum displacement (in the $x$-direction) and the relative error in terms of the absolute difference between the simulated value from the reference value ($\Bar{u}_x^{\textrm{ref}}=0.0625$[m]). In analyzing the impact of the child mesh resolution, the parent resolution is fixed and the relative resolution ($r^\textrm{p-c}$) is varied. Note that the relative resolution is assumed $r^{\textrm{p-c}}\sim\mathcal{O}(1)$ to prevent erroneous conclusions resulting either due to computational instabilities (when $r^\textrm{p-c} \ll 1$) or due to poor efficiency (when $r^\textrm{p-c} \gg 1$). The full-field simulation for a specific combination of the two kernels, nonlocal horizon length, and mesh sizes, are presented in Fig.~(\ref{fig: verification_case_1}) and Fig.~(\ref{fig: verification_case_2}). The excellent match between the simulated and reference responses successfully validates the implementation of M${}^2$-FEM (particularly, the mesh-decoupling and scale-bridging schemes, that are the pillars of M${}^2$-FEM). 

\begin{table}[ht!]
\caption{Relative error in the maximum displacement along $x$-direction obtained by using $\text{M}^2$-FEM, against the reference value, parameterized for different values of the nonlocal horizon length ($l^{nl}$) and kernel decay rate ($\tau$) for \textit{Case study 1}. The effect of the (a) parent domain resolution and (b) relative resolution on accuracy are also analyzed.}
\label{tab: Exp_accuracy}
    \centering
    \begin{tabular}{c|c|c c|c c|c|c c|c c}
    \hline\hline
    \multirow{2}{*}{$l^\textrm{nl}$}
    & \multirow{2}{*}{(a) $M^\textrm{p}$}
    & \multicolumn{2}{c}{$\tau = 0.002L$}\vline
    & \multicolumn{2}{c}{$\tau = 0.003L$}\vline
    & \multirow{2}{*}{(b) $r^\textrm{p-c}$} 
    & \multicolumn{2}{c}{$\tau = 0.002L$}\vline
    & \multicolumn{2}{c}{$\tau = 0.003L$} \\
    \cline{3-6}
    \cline{8-11}
    & & $\overline{u}_x\times100$ & Er[\%]
    & $\overline{u}_x\times100$ & Er[\%]
    & & $\overline{u}_x\times100$ & Er[\%]
    & $\overline{u}_x\times100$ & Er[\%]\\
    \hline
    \multirow{3}{*}{$l^\textrm{nl}=0.2L$}
    & $8$ & 5.96 & 4.67 & 6.0647 & 2.96 & $1$ & 6.23 & 0.40 & 6.24 & 0.24 \\
    & $16$ & 6.23 & 0.28 & 6.2450 & 0.08 & $2$ & 6.24 & 0.17 & 6.24 & 0.23 \\
    & $24$ & 6.24 & 0.11 & 6.2377 & 0.12 & $3$ & 6.24 & 0.17 & 6.24 & 0.23 \\
    \hline
    \multirow{3}{*}{$l^\textrm{nl}=0.3L$}
    & $8$ & 5.96 & 4.65 & 6.087 & 2.64 & $1$ & 6.23 & 0.40 & 6.24 & 0.24 \\
    & $16$ & 6.23 & 0.28 & 6.24 & 0.08 & $2$ & 6.24 & 0.17 & 6.24 & 0.23 \\
    & $24$ & 6.24 & 0.11 & 6.4 & 0.20 & $3$ & 6.24 & 0.17 & 6.24 & 0.23 \\
    \hline\hline
    \end{tabular}
\end{table}

\begin{table}[ht!]
\caption{Relative error in the maximum displacement along $x$-direction obtained by using $\text{M}^2$-FEM, against the reference value, parameterized for different values of the nonlocal horizon length ($l^{\textrm{nl}}$) and power-law orders ($\alpha$) for \textit{Case study 2}. The effect of the (a) parent domain resolution and (b) relative resolution on accuracy are also analyzed.}
\label{tab: PL_accuracy}
    \centering
    \begin{tabular}{c|c|c c|c c|c|c c|c c}
    \hline\hline
    \multirow{2}{*}{$l^\textrm{nl}$}
    & \multirow{2}{*}{(a) $M^\textrm{p}$}
    & \multicolumn{2}{c}{$\alpha=0.6$}\vline
    & \multicolumn{2}{c}{$\alpha=0.3$}\vline
    & \multirow{2}{*}{(b) $r^\textrm{p-c}$} 
    & \multicolumn{2}{c}{$\alpha=0.6$}\vline
    & \multicolumn{2}{c}{$\alpha=0.3$} \\
    \cline{3-6}
    \cline{8-11}
    & & $\overline{u}_x\times100$ & Er[\%]
    & $\overline{u}_x\times100$ & Er[\%]
    & & $\overline{u}_x\times100$ & Er[\%]
    & $\overline{u}_x\times100$ & Er[\%]\\
    \hline
    \multirow{3}{*}{$l^\textrm{nl}=0.3L$}
    & $8$ & 6.15 & 1.52 & 6.24 & 0.15 & $1$ & 6.26 & 0.19 & 6.22 & 0.43 \\
    & $16$ & 6.20 & 0.75 & 6.28 & 0.45 & $2$ & 6.22 & 0.34 & 6.21 & 0.68 \\
    & $24$ & 6.30 & 0.80 & 6.13 & 1.87 & $3$ & 6.223 & 0.43 & 6.21 & 0.70 \\
    \hline
    \multirow{3}{*}{$l^\textrm{nl}=0.5L$}
    & $8$ & 6.34 & 1.51 & 6.48 & 3.60 & $1$ & 6.28 & 0.42 & 6.35 & 1.64 \\
    & $16$ & 6.27 & 0.40 & 6.32 & 1.19 & $2$ & 6.23 & 0.24 & 6.27 & 0.34 \\
    & $24$ & 6.26 & 0.18 & 6.29 & 0.59 & $3$ & 6.23 & 0.36 & 6.25 & 0.08 \\
    \hline\hline
    \end{tabular}
\end{table}

\subsection{Convergence analysis}
\label{sec: 4.3}
Analogous to existing integral solvers \cite{pisano2009nonlocal,sidhardh2018effect,patnaik2020ritz,patnaik2021fractional}, the convergence of the $\text{M}^2$-FEM involves a more elaborate analysis (when compared to classical FEM) due to the additional child domain integral approximation that is nested within the weak-form integral. In addition to completeness and compatibility conditions typical of the finite element simulation of local elasticity problems \cite{oden2012introduction}, the convergence of the nested parent-child integrals in the simulation of integral elasticity theories is affected also by the relative resolution of the parent and child meshes (equivalent to the dynamic rate of convergence in mesh-coupled solvers \cite{sidhardh2018effect,patnaik2020ritz}; shown below) and the nature of the kernel function (more particularly, the spatial decay rate and its resolution within the child mesh) \cite{patnaik2022displacement}. We present comprehensive parametric studies by leveraging both $\textit{h}$-refinement and the $\textit{p}$-refinement to investigate the convergence conditions in terms of the aforementioned factors. In the subsequent analysis, we define the \textit{convergence criterion} as a (arbitrary and commonly adopted) threshold of 2\% relative difference in the Euclidean distance between the displacement fields simulated from successive spatial refinements (in either the parent or child domain); the relative difference is denoted as $\Delta$ simply, without any super- or sub-scripts (to avoid confusion with the parent and child element sizes).

\paragraph*{$\textit{\textbf{h}}$-refinement}\mbox{}\\
We perform $\textit{h}$-refinement in the child domain because the analysis of accuracy presented in \S\ref{sec: 4.2}, as well as convergence studies conducted for existing integral solvers \cite{patnaik2020ritz,patnaik2021fractional}, indicate that the performance of the integral solver is sensitive to the child domain resolution. In fact, we first show that the 2D child domain resolution is equivalent to the dynamic rate of convergence that explicitly determines the convergence in mesh-coupled solvers:
\begin{equation}
    \left[\frac{l^{\mathrm{nl}}}{\Delta^\mathrm{c}_x}\right] \times \left[\frac{l^{\mathrm{nl}}}{\Delta^\mathrm{c}_y}\right] \equiv \frac{1}{r^{\textrm{p-c}}} \left[ \frac{l^{\mathrm{nl}}}{\Delta^\mathrm{p}_x} \right] \times \frac{1}{r^{\textrm{p-c}}} \left[\frac{l^{\mathrm{nl}}}{\Delta^\mathrm{p}_y} \right] \xRightarrow{r^{\textrm{p-c}}=1} \underbrace{\left[ \frac{l^{\mathrm{nl}}}{\Delta^\mathrm{p}_x} \right]}_{\mathcal{N}_x} \times \underbrace{  \left[\frac{l^{\mathrm{nl}}}{\Delta^\mathrm{p}_y} \right]}_{\mathcal{N}_y}
\end{equation}
In the second equivalence above, we assume $r^{\textrm{p-c}}=1$, as is the case for mesh-coupled solvers, and $\mathcal{N}_x \times \mathcal{N}_y$ is the dynamic rate of convergence. Remarkably, this mathematical deduction suggests that the overall convergence (in addition to the desired accuracy) can be directly achieved by independently tuning the child mesh resolution while still maintaining a coarse parent mesh for superior efficiency (thanks, to the decoupling; see \S\ref{sec: 5.2}).

The results of a parametric convergence analysis corresponding to the case studies introduced previously is presented in Table~\ref{tab: convergence_1}. As evident from the results, the convergence for the exponential kernel (case study 1) is achieved at $\mathcal{N}_x\times\mathcal{N}_y=6\times6$ and for the power-law kernel at $\mathcal{N}_x\times\mathcal{N}_y=5\times5$ (case study 2). In other terms, atleast 6 and 5 elements are required (on an average) along each direction to achieve a detailed resolution of the displacement profile following the exponential and power-law kernels, respectively. The (slightly) slower convergence observed in the exponential kernel, when compared to the power-law kernel, on refinement is due to the faster spatial decay of the exponential kernel which requires a more detailed mesh for better resolution; see similar analysis also in \cite{patnaik2022displacement}. However, in comparison to existing solvers which attain convergence around $\mathcal{N}_{\square}\sim\mathcal{O}(10)$ (where, $\square\in\{x,y\}$) \cite{patnaik2021fractional,pisano2009nonlocal,patnaik2020ritz,sidhardh2018effect}, $\text{M}^2$-FEM rapid attains convergence on successive mesh refinement. The latter remarkable outcome is a direct impact of the powerful ability of the decoupled child mesh in enforcing a convergence (from the top, as in local FEM \cite{oden2012introduction}) to the real solution independently of the parent mesh; due to the minimal transmission of approximation error to and from the parent mesh.

\begin{table}[ht!]
\caption{Convergence of the $\text{M}^2$-FEM simulation of the nonlocal problems in the two case studies via $\textit{h}$-refinement.}
    \label{tab: convergence_1}
    \centering
    \resizebox{\textwidth}{!}{
    \begin{tabular}{c|c|c c|c c|c|c|c c|c c}
    \hline\hline
    \multicolumn{6}{c}{\textit{Case study 1}}\vline
    & \multicolumn{6}{c}{\textit{Case study 2}} \\
    \hline
    \hline
    \multirow{2}{*}{$l^\textrm{nl}$}
    & \multirow{2}{*}{$\mathcal{N}_x\times\mathcal{N}_y$}
    & \multicolumn{2}{c}{$\tau=0.002L$}\vline
    & \multicolumn{2}{c}{$\tau=0.003L$}\vline
    & \multirow{2}{*}{$l^\textrm{nl}$}
    & \multirow{2}{*}{$\mathcal{N}_x\times\mathcal{N}_y$}
    & \multicolumn{2}{c}{$\alpha=0.6$}\vline
    & \multicolumn{2}{c}{$\alpha=0.3$} \\
    \cline{3-6}
    \cline{9-12}
    & & $\overline{u}_x\times100$ & $\Delta$ [\%]
    & $\overline{u}_x\times100$ & $\Delta$ [\%]
    & & & $\overline{u}_x\times100$ & $\Delta$ [\%]
    & $\overline{u}_x\times100$ & $\Delta$ [\%] \\
    \hline
    \multirow{4}{*}{$0.2L$}
    & $2\times2$ & 6.12 & - & 6.13 & - 
    & \multirow{4}{*}{$0.25L$}
    & $2\times2$ & 6.44 & - & 6.40 & - \\
    & $4\times4$ & 6.27 & 2.59 & 6.27 & 2.32 & & $4\times4$ & 6.30 & 2.19 & 6.30 & 1.56 \\
    & $5\times5$ & 6.24 & 0.46 & 6.24 & 0.45 & & $5\times5$ & 6.29 & 0.15 & 6.29 & 0.20 \\
    & $6\times6$ & 6.25 & 0.24 & 6.25 & 0.25 & & $6\times6$ & 6.28 & 0.09 & 6.28 & 0.12 \\
    \hline
    \multirow{4}{*}{$0.4L$}
    & $2\times2$ & 3.26 & - & 4.32 & -
    & \multirow{4}{*}{$0.5L$}
    & $2\times2$ & 6.76 & - & 6.92 & - \\
    & $4\times4$ & 6.12 & 87.62 & 6.13 & 42.13 & & $4\times4$ & 6.39 & 5.54 & 6.42 & 7.27 \\
    & $6\times6$ & 6.22 & 1.79 & 6.23 & 1.72 & & $5\times5$ & 6.35 & 0.63 & 6.37 & 0.86 \\
    & $8\times8$ & 6.27 & 0.78 & 6.27 & 0.59 & & $6\times6$ & 6.32 & 0.35 & 6.34 & 0.47 \\
    \hline\hline
    \end{tabular}}
\end{table}

\paragraph*{$\textit{\textbf{p}}$-refinement}\mbox{}\\
For a complete analysis of the convergence conditions derived from $\textit{h}$-refinement, we additionally perform a $\textit{p}$-refinement fixing the parent mesh parameters at $M^\textrm{p}=12$ and $r^\textrm{p-c} = 1$. For this purpose, $1^\textrm{st}$ and $2^\textrm{nd}$ order Lagrange shape functions are chosen. Note that, here the $\textit{p}$-refinement is only applicable to the parent domain since the DOF (and hence the shape functions), following Eringen's approach are defined at the parent domain. The results of this analysis are presented in Table \ref{tab:convergence_2}. The nonlocal parameters assumed here are different from that in $\textit{h}$-refinement, to verify the universality of the convergence conditions. As evident from Table \ref{tab:convergence_2}, the relative change in the Euclidean distance between the displacement fields obtained using the $1^\textrm{st}$ and $2^\textrm{nd}$ order Lagrange shape functions is less than the \textit{a priori} defined threshold (that is, $\Delta < 2\%$), indicating converged simulations.

\begin{table}[ht!]
\caption{Convergence of $\text{M}^2$-FEM solution to the nonlocal problem in the two case studies using $\textit{p}$-refinement.}
    \label{tab:convergence_2}
    \centering
    \begin{tabular}{c|c|c c c|c|c|c c c}
    \hline\hline
    \multicolumn{5}{c}{\textit{Case study 1}}\vline
    & \multicolumn{5}{c}{\textit{Case study 2}} \\
    \hline
    \hline
    \multirow{2}{*}{$l^\textrm{nl}$} & \multirow{2}{*}{$\tau$} & $1^\textrm{st}$ order & $2^\textrm{nd}$ order & $\Delta$ & \multirow{2}{*}{$l^\textrm{nl}$} & \multirow{2}{*}{$\alpha$} & $1^\textrm{st}$ order & $2^\textrm{nd}$ order & $\Delta$ \\
    & & $\overline{u}_x\times100$ & $\overline{u}_x\times100$ & [\%] & & & $\overline{u}_x\times100$ & $\overline{u}_x\times100$ & [\%] \\
    \hline
    \multirow{2}{*}{$0.2L$} & $0.002L$ & 6.19 & 6.15 & 0.65 & \multirow{2}{*}{$0.3L$} & $0.6$ & 6.34 & 6.27 & 1.10 \\
    & $0.003L$ & 6.26 & 6.28 & 0.32 & & $0.3$ & 6.37 & 6.18 & 2.98 \\
    \hline
    \multirow{2}{*}{$0.3L$} & $0.002L$ & 6.08 & 6.15 & 1.15 & \multirow{2}{*}{$0.5L$} & $0.6$ & 6.36 & 6.28 & 1.26 \\
    & $0.003L$ & 6.21 & 6.29 & 1.29 & & $0.3$ & 6.42 & 6.35 & 1.09 \\
    \hline\hline
    \end{tabular}
\end{table}

\clearpage
\section{Impact of \texorpdfstring{$\text{M}^2$}--FEM on multiscale modeling}
\label{sec: 5}
In order to highlight the significant impact of M${}^2$-FEM on multiscale modeling, we showcase the ability of the underlying mesh-decoupling technique to: 1) enable the use of unstructured meshes for the simulation of non-standard non-rectangular geometries with (geometrically) non-similar nonlocal horizons; and 2) achieve high accuracy and high efficiency simultaneously, by tuning the decoupled child and parent mesh-sizes, respectively.

\subsection{Simulation of non-rectangular multiscale configurations}
\label{sec: 5.1}
\begin{figure}[ht!]
    \centering
    \includegraphics[width=1\linewidth]{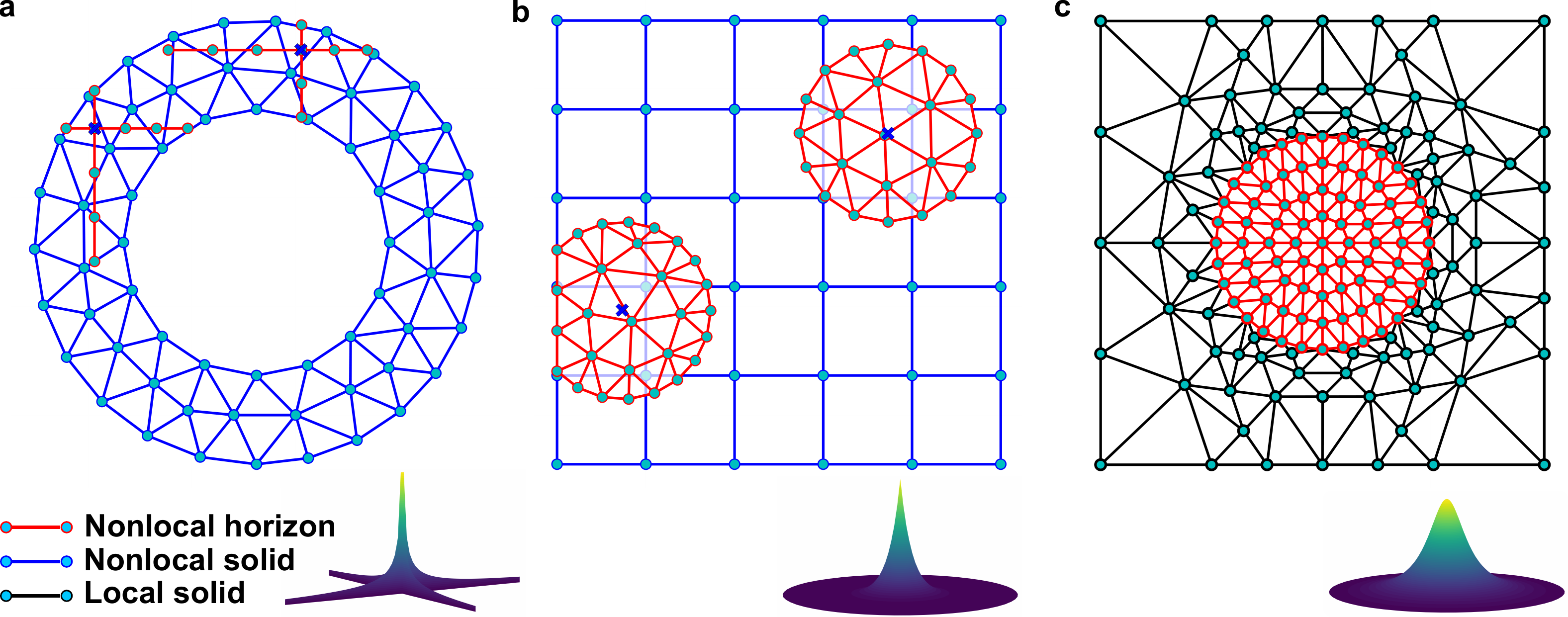}
    \caption{Schematic illustrating mesh and nonlocal kernel function configurations in the three numerical test cases. \textbf{a} Nonlocal annulus: the annular parent domain is discretized by unstructured triangular elements while the child domain consisting of two orthogonal lines is discretized by 1D line elements. \textbf{b} Nonlocal chip: the rectangular parent domain and the circular child domain are discretized by structured quadrilateral elements and unstructured triangular elements, respectively. \textbf{c} Nonlocal inclusion: unstructured triangular mesh is used to discretize the whole domain and adaptive mesh refinement technique is used in the nonlocal inclusion to achieve higher accuracy. Different nonlocal kernel functions applied in the three simulation cases are presented at the bottom-right corner of each configuration.}
    \label{fig: unstructured_mesh}
\end{figure}

We consider a set of solids with non-rectangular shape and non-rectangular nonlocal horizons. These choices enforce the use of unstructured mesh in both the parent and child domains and hence, (recall from \S\ref{sec: introduction}) cannot be simulated by existing integral solvers. Three configurations, illustrated in Fig.~(\ref{fig: unstructured_mesh}) along with their unstructured mesh, are considered. The three configurations were chosen following a careful survey of some common geometrical features underpinning multiscale solids in practical applications; they are not defined arbitrarily to merely show the need for unstructured meshing. Further, the nonlocal kernels used to simulate the configurations are also very closely representative (in a functional sense) of the actual microstructural signature in the emerging nonlocal effects noted in corresponding applications. For all the configurations, the nonlocal horizon is truncated at points lying close to and on the material boundaries \cite{polizzotto2001nonlocal,patnaik2020generalized} resulting in non-constant irregular shapes of the horizon. This phenomenon is schematically illustrated for selected points in Fig.~(\ref{fig: unstructured_mesh}).

\subsubsection{Relevance and behavior of the selected configurations}
The geometry, practical relevance, and mechanical behavior of the three configurations are discussed below: 
\begin{itemize}[leftmargin=*]
    \item \textit{Nonlocal annulus}: a nonlocal annular domain at the parent mesh level. At each point within the annulus, the nonlocal horizon is assumed to be restricted along bidirectional line segments. From a numerical perspective, the differential dimensions assumed for the parent (2D) and child (1D) highlight the significant flexibility afforded by $\text{M}^2$-FEM; indeed this choice presents the most drastic case of geometric non-similarity between the parent and child domains. Such reduced-order and preferential directions of localization of nonlocal effects is typical of orthotropic and anisotropic microstructures in composites and porous solids \cite{kanoute2009multiscale,mannan2018stiffness}. The overall nonlocal annulus geometry is commonly noted in naturally-occurring porous solids (such as, animal bones and bamboo) and several real-world applications including energy conversion devices (batteries) and biomedical implants \cite{mannan2018stiffness,mannan2017correlations}. As shown in \cite{sebaa2006application,patnaik2022role,patnaik2020generalized}, power-law kernels closely capture the deformation behavior of these solids.

    \item \textit{Nonlocal plane strain domain}: a square-shaped parent domain whose material points are associated with a circular nonlocal horizon (geometrical non-similar to the square-shaped parent domain). This configuration could be very relevant to nano- and micro- electromechanical systems \cite{sidhardh2018effect}, semiconductor devices, and even periodic metastructures embedded with acoustic black holes \cite{nair2021nonlocal}. Exponential kernel functions are more representative of the dispersive behavior in these structural configurations \cite{sidhardh2018effect,nair2021nonlocal}.
    
    \item \textit{Local plane strain domain with nonlocal inclusion}: a square-shaped plane strain local domain with an embedded circular nonlocal inclusion. The inclusion is fully nonlocal, meaning that the nonlocal horizon at any material point within the inclusion, coincides with the entire inclusion. This configuration is an excellent example of interface problems and can be seen as a multiscale extension of the classical inclusion problem and can be found commonly in functionally-graded solids, layered media (such as, soil), and solids with periodic reinforcements, interstitial defects, or distributed damage \cite{sidhardh2018inclusion,romkes2004adaptive}. The strong nonlocality and dispersion noted in these applications is well captured by rational kernels \cite{sidhardh2018inclusion,seleson2013interface}.
\end{itemize}

\begin{figure}[htbp]
    \centering
    \includegraphics[width=\linewidth]{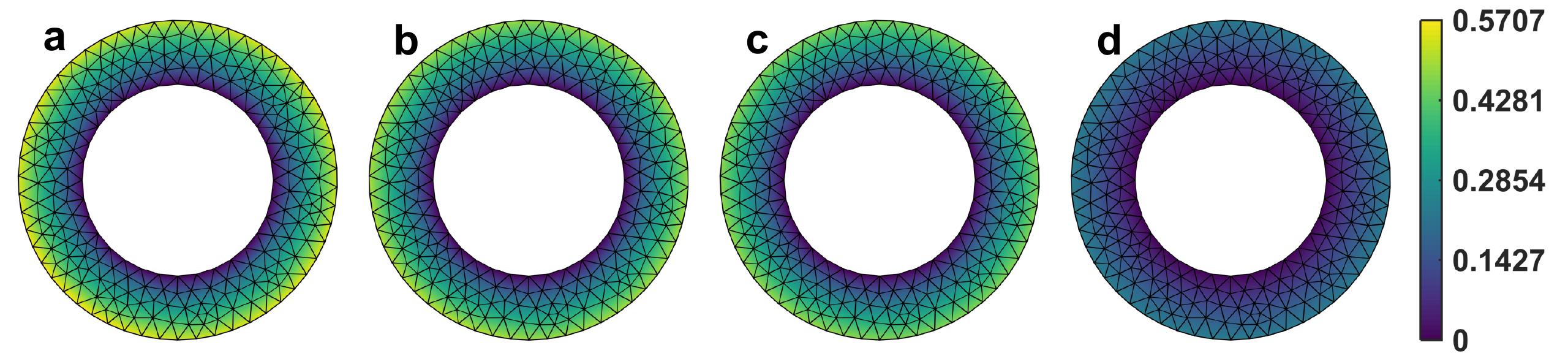}
    \caption{Response of the nonlocal annulus in terms of the absolute displacement field, simulated via $\text{M}^2$-FEM. The results are parameterized for different values of the nonlocal strength parameter $\alpha$: \textsb{a} $\alpha=0.4$, \textsb{b} $\alpha=0.6$, \textsb{c} $\alpha=0.8$, and \textsb{d} $\alpha=1$.}
    \label{fig: solution_case_1}
\end{figure}

\paragraph*{Nonlocal annulus.}
The inner and outer radii of the annulus are taken as $r^{\textrm{in}}=0.3\textrm{m}$ and $r^{\textrm{out}}=0.5\textrm{m}$, respectively. Further, the length scale of the nonlocal horizon is taken as $l^\textrm{nl}=0.2\textrm{m}$. The parent and child domains are discretized using triangular and linear Lagrange elements, respectively. The kernel function is defined (similar to \cite{patnaik2022role}) as: 
\begin{subequations}\label{eq: kernel_case_1}
\begin{align}
    \mathcal{K}(x,y,x^\prime,y^\prime,\alpha) &= \frac{1}{2}\left[\mathcal{K}_x(x,y,x^\prime,y^\prime,\alpha) + \mathcal{K}_y(x,y,x^\prime,y^\prime,\alpha)\right] \\
    \mathcal{K}_x(x,y,x^\prime,y^\prime,\alpha) &=
    \frac{1}{2\Gamma(1-\alpha)}\frac{1}{\vert x-x^\prime \vert^\alpha}, \quad \vert x-x^\prime \vert < l^\textrm{nl},~y = y^\prime \\
    \mathcal{K}_y(x,y,x^\prime,y^\prime,\alpha) &=  \frac{1}{2\Gamma(1-\alpha)}\frac{1}{\vert y-y^\prime \vert^\alpha}, \quad \vert y-y^\prime \vert < l^\textrm{nl},~x = x^\prime
\end{align}
\end{subequations}
where $\mathcal{K}_x$ and $\mathcal{K}_y$ are singular power-law functions defined in the $x$ and $y$ directions. The above definition reduces the 2D nonlocal horizon to bidirectional 1D line segments. The inner boundary of the annulus is fixed, and the outer boundary is loaded with normal tractions of constant magnitude $\vert \bm{t}^{\textrm{out}} \vert = 0.1\textrm{Pa}$. The response of the nonlocal annulus, in terms of the absolute value of the displacement field, parameterized for different values of the fractional-order $\alpha$, are presented in Fig.~(\ref{fig: solution_case_1}).
As evident, $\text{M}^2$-FEM enables the smooth (lack of spurious gradients) and non-divergent simulation of non-rectangular geometries. Observe that the magnitude of the maximum absolute displacement (around the outer boundary) increases as the power-law order $\alpha$ decreases. This simulated softening behavior is consistent with the predictions of fractional-order approaches to nonlocal elasticity \cite{ding2022multiscale}. Note that the power-law kernel in Eq.~\eqref{eq: kernel_case_1} is denominated constant-order kernel in FC \cite{das2011functional}. To this end, we highlight that $\text{M}^2$-FEM can also be used to simulate variable- and distributed-order kernels that lead to functionally-varying and functionally-distributed space of the power-law exponent, respectively. These advanced kernels were used recently to develop the highly generalized variable-order and distributed-order nonlocal approaches \cite{patnaik2021variable,ding2022multiscale}. The asymmetric and anisotropic kernels (unlike, Eringen's approach \cite{batra159misuse}) used in these approaches enables the simulation strongly nonlocal and multiscale effects induced by intricate microstructures \cite{patnaik2022distillation,ding2022multiscale2}.
 
\begin{figure}[ht!]
    \centering
    \includegraphics[width=\linewidth]{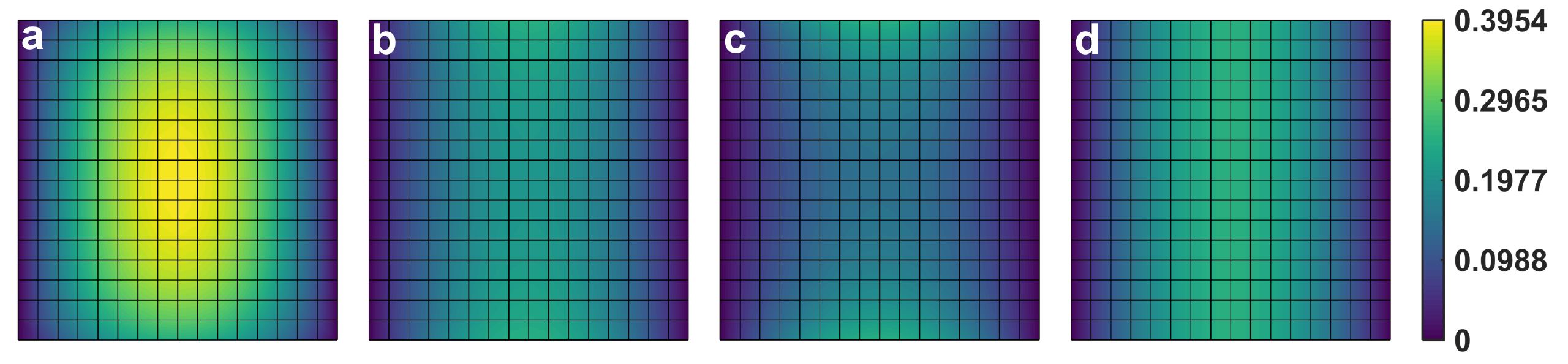}
    \caption{Response of the nonlocal plane strain domain simulated via $\text{M}^2$-FEM, is presented in terms of the displacement fields along $x$ direction ($u_x$). The response is parameterized for different values of the nonlocal strength $\tau_2$: \textsb{a} $\tau_1=1/1000$, \textsb{b} $\tau_1=1/2000$, and \textsb{c} $\tau_1=1/3000$, with a fixed value of $\tau_2=1/100$. The local response is presented in $\textsb{d}$ for reference.}
    \label{fig: solution_case_2_tau_1}
\end{figure}

\begin{figure}[ht!]
    \centering
    \includegraphics[width=\linewidth]{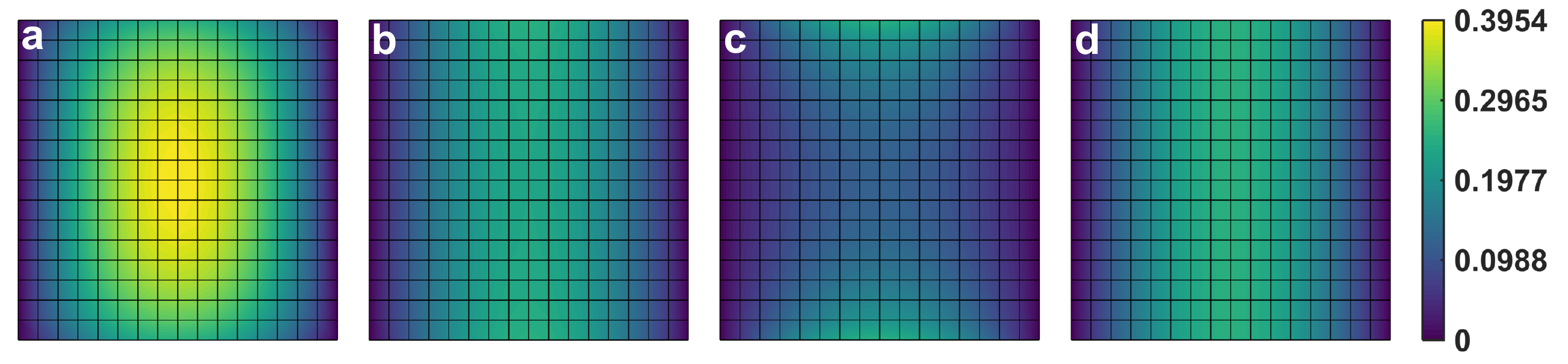}
    \caption{Response of the nonlocal plane strain domain simulated via $\text{M}^2$-FEM, is presented in terms of the displacement fields along $x$ direction ($u_x$). The response is parameterized for different values of nonlocal strength parameter $\tau_2$: \textsb{a} $\tau_2=1/100$, \textsb{b} $\tau_2=1/75$, and \textsb{c} $\tau_2=1/50$, with a fixed value of $\tau_1=1/1000$. The local response is presented in $\textsb{d}$ for reference.}
    \label{fig: solution_case_2_tau_2}
\end{figure}

\paragraph*{Nonlocal plane strain domain.}
This configuration, as shown in Fig.~(\ref{fig: unstructured_mesh}-b), consists of a square-shaped parent domain with size $L=1\textrm{m}$, and a circular child domain with radius $l^\textrm{nl} = 0.2\textrm{m}$. The parent domain is discretized uniformly using structured quadrilateral elements and the child domain is discretized randomly using unstructured triangular elements. Note that the truncation of nonlocal horizon near boundaries, requires the regeneration and re-descretization of the corresponding child mesh. We adopt a 2D and radially symmetric exponential kernel:
\begin{equation}\label{eq: kernel_case_2}
    \mathcal{K}(\bm{x},\bm{x^\prime}) = \frac{1}{\tau_1}\mathrm{exp}\left(\frac{\vert\bm{x}-\bm{x}^\prime\vert}{\tau_2}\right), \quad \vert\bm{x}-\bm{x}^\prime\vert < l^\textrm{nl}
\end{equation}
where $\tau_1$ and $\tau_2$ are constitutive parameters that govern the amplitude and decay rate of the kernel, respectively. Results of two parametric studies within this example, for different values of $\tau_1$ and $\tau_2$, are presented in Fig.~(\ref{fig: solution_case_2_tau_1}) and Fig.~(\ref{fig: solution_case_2_tau_2}), respectively. Aside from the present discussion, an interesting result emerges from the two studies. The results in Figs.~(\ref{fig: solution_case_2_tau_1},\ref{fig: solution_case_2_tau_2}) show that, depending on the value of parameters $\tau_1$ and $\tau_2$, the nonlocal kernel in Eq.~(\ref{eq: kernel_case_2}), can represent either structural softening or stiffening behavior (a feature often difficult to implement in multiscale constitutive modeling\footnote{A discussion on the constitutive modeling of stiffening-softening effects is beyond the scope of the present study. The interested reader is referred to \cite{eringen1984theory,polizzotto2003unified,patnaik2022displacement,rahimi2020non} for detailed discussions on the ability of nonlocal constitutive theories in modeling softening and/or stiffening effects.}). More specifically, results in Fig.~(\ref{fig: solution_case_2_tau_1}-a) and Fig.~(\ref{fig: solution_case_2_tau_2}-a) show the emergence of the softening behavior. The softening effect is evident from the larger displacement amplitude (around the central region of the domain) when compared with the local counterpart in Fig.~(\ref{fig: solution_case_2_tau_1}-d) and Fig.~(\ref{fig: solution_case_2_tau_2}-d). On the other hand, a stiffening behavior results in smaller displacement amplitudes (compared to the local response) as noted in Figs.~(\ref{fig: solution_case_2_tau_1}-c,\ref{fig: solution_case_2_tau_2}-c).

\begin{figure}[ht!]
    \centering
    \includegraphics[width=\linewidth]{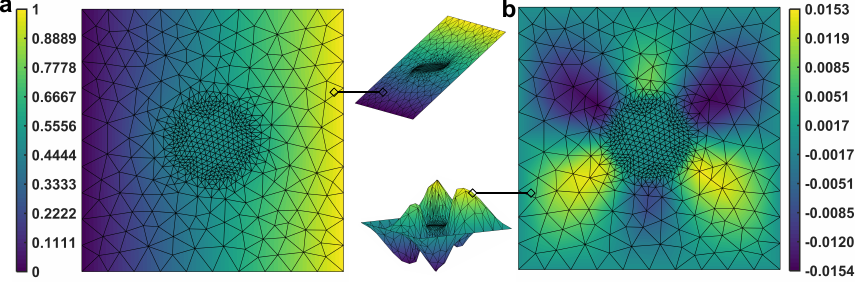}
    \caption{Simulation of the local plane strain domain with nonlocal inclusion and nonlocal interface, via $\text{M}^2$-FEM. Results are presented in terms of the displacement field along the: $\textsb{a}$ $x$-direction ($u_x$), and $\textsb{b}$ y-direction ($u_y$). Orthographic views of the displacement fields are also presented to enable a complete visualization.}
    \label{fig: solution_case_3}
\end{figure}

\paragraph*{Local plane strain domain with nonlocal inclusion.}
The outer dimensions of the local matrix and the radius of the inner inclusion are assumed as $L_x \times L_y = 1\textrm{m} \times 1\textrm{m}$ and $r^\textrm{in} = 0.15\textrm{m}$, respectively. Unstructured meshes with triangular elements were used to discretize the matrix and the inclusion (including both parent and child domains). Notably, mesh-decoupling directly enabled the use of mesh-grading and mesh-refinement (see Fig.~(\ref{fig: solution_case_3})) \cite{oden2012introduction}, to reduce the overall computational cost. More specifically, these techniques reduced the computational resources dedicated to the local matrix, while simultaneously ensuring a detailed resolution of the nonlocal inclusion and the interface. This result is of practical relevance in the simulation of more complex configurations, such as those with distributed inclusions. Now, to simulate this configuration we assume a rational kernel:
\begin{equation}
    \mathcal{K}(\bm{x},\bm{x^\prime}) = \frac{1}{\tau_1 \vert\bm{x}-\bm{x}^\prime\vert^2+\tau_2}, \quad \vert\bm{x}-\bm{x}^\prime\vert \leq 2r^\textrm{in}
\end{equation}
where $\tau_1=100$ and $\tau_2=5000$ are nonlocal model parameters. Displacement boundary conditions are employed:
\begin{equation}
    \bm{u}(x,y) \equiv \{u_x(x,y),u_y(x,y)\} = \{x,0\} \quad \forall (x,y)
\end{equation}
The simulated response is presented Fig.~(\ref{fig: solution_case_3}). Notice that the presence of nonlocal effects distorts the displacement field within an appreciable neighbourhood surrounding the interface. A clear change in displacement gradient, also called 'weak discontinuity', is observed along the circumference of the interface. This behavior is typical of the response of nonlocal inclusions \cite{patnaik2022role,sidhardh2018inclusion}, and it is well captured by $\text{M}^2$-FEM. The results are remarkable because, unlike existing approaches, $\text{M}^2$-FEM readily discretizes the circular interface between the nonlocal inclusion and the local matrix independently, while simultaneously capturing the multiscale interaction via scale-bridging.

\subsection{High accuracy with high efficiency: overcoming restrictive trade-offs}
\label{sec: 5.2}
The most significant impact of $\text{M}^2$-FEM, in view of the highly multiscale character of emerging \cite{patnaik2022distillation,failla2020advanced} and future applications \cite{van2020roadmap,hoekstra2014multiscale}, is the ability to simultaneously achieve high accuracy and high efficiency while requesting limited computational resource. This is in sharp contrast to traditional solvers that impose implicit trade-offs between accuracy and efficiency. Note that, differently from the approximation of local elastic potential energies, multiscale constitutive models require the approximation of additional integrals at each Gauss point in the parent domain, which drastically increases time complexity of these nonlocal problems\footnote{Note that the same concept also extends to gradient theories where additional gradients would necessitate numerical evaluation. Although beyond the scope of this study, recall from \S\ref{sec: introduction} that, gradient models are less general when compared to integral models. In regards to the present discussion, they also lead to practical difficulties in the numerical implementation, rooted in the continuity of evaluated higher-order gradients.}. As explained in detail in \S\ref{sec: introduction}, sharing the mesh between the parent and child domains drastically increases the overall computational cost when tuning the element size mesh for better accuracy in the nonlocal integrals\footnote{As seen in \S\ref{sec: 4.2} and literature \cite{patnaik2020ritz,pisano2009nonlocal}, the numerical accuracy of nonlocal finite element approaches is more sensitive to the approximation of nonlocal integrals associated with the fine material scales.}. On the other hand, the decoupling of the meshes in $\text{M}^2$-FEM, allows an independent tuning of the child element size for better accuracy (via better approximation of the fine-scale nonlocal integrals) and the parent element size for better efficiency. To enable a more transparent visualization of this concept, we analyze the time complexity of both $\text{M}^2$-FEM and traditional nonlocal finite element approaches when simulating integral elasticity problems.

\begin{figure}[ht!]
    \centering
    \includegraphics[width=\linewidth]{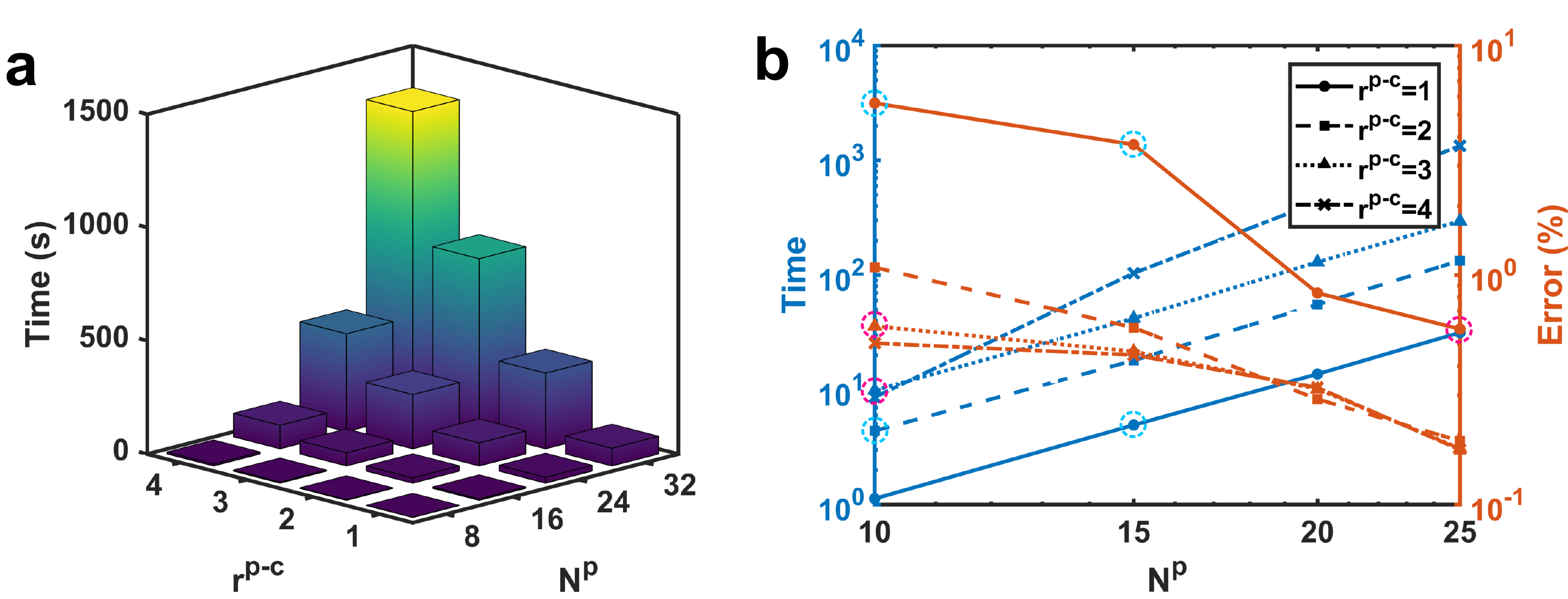}
    \caption{Variation of both computational time and error with respect to the parent and child mesh numbers. \textbf{a} 3D bar plot showing growth in computational time following an increase in both the parent and child mesh numbers. \textbf{b} Variation of computational time cost (left $y$-axis and blue curves) and the relative error (right $y$-axis and orange curves) with different mesh numbers. The legend specifying different line styles applies to both blue and orange curves.}
    \label{fig: computational_analysis}
\end{figure}

We consider a benchmark problem similar to that shown in Fig.~(\ref{fig: FE_configuration}). More specifically, we consider a square parent domain (size $L \times L$) which is discretized uniformly by the parent mesh with mesh number $M^\textrm{p}$ in both the $x$ and $y$ directions. Further, the child domain at each parent Gauss point is also assumed to be square-shaped (size $l^\textrm{nl} \times l^\textrm{nl}$), and it is discretized uniformly by the child mesh with mesh number $M^\textrm{c}$ in both the $x$ and $y$ directions. To simplify the overall analysis, we do not consider the truncation of the child domain caused by the presence of the material boundary and hence, ignore the change of child mesh number for the parent Gauss points lying close to the boundaries. For the above defined configuration, we obtain the total number of parent Gauss points in the parent domain and the total number of child Gauss points in each child domain, as:
\begin{subequations}
\begin{align}
    M^{g^\textrm{p}} &= g^\textrm{p}\left(M^\textrm{p}\right)^2 \\
    M^{g^\textrm{c}} &= g^\textrm{c}\left(M^\textrm{c}\right)^2
\end{align}
\end{subequations}
where $g^\textrm{p}$ and $g^\textrm{c}$ are the number of parent and child Gauss points in each parent and child element, respectively. Now, the total number of individual evaluations and summations performed in the approximation of the nonlocal integration and the weak-form integration, both following the numerical quadrature rules, can be approximated (accurate to the limit of multiplicative or additive constants) as:
\begin{equation}
    M^{\textrm{sum}} \approx g^\textrm{p}g^\textrm{c}\left(M^\textrm{c}\right)^2\left(M^\textrm{p}\right)^2
\end{equation}
Note that the above expression is valid for all existing nonlocal finite element approaches (including, $\text{M}^2$-FEM) for a first-order (one) parent-child nesting; see similar analysis in \cite{patnaik2021fractional}. With the above understanding, recall the relative resolution $r^\textrm{p-c}$ of $\text{M}^2$-FEM, used in \S\ref{sec: 4.1} to represent the ratio between the parent and child mesh size. By using the definitions of the relative resolution, and the parent and child mesh numbers, $M^{\textrm{sum}}$ can be expressed as:
\begin{equation}\label{eq: M_sum}
    M^{\textrm{sum}} 
    \approx g^\textrm{p}g^\textrm{c}\left(r^\textrm{p-c}\right)^2 \left(\frac{l^\textrm{nl}}{L}\right)^2 \left(M^\textrm{p}\right)^4
    \equiv g^\textrm{p}g^\textrm{c} \frac{1}{\left(r^\textrm{p-c}\right)^2} \left(\frac{L}{l^\textrm{nl}}\right)^2 \left(M^\textrm{c}\right)^4
\end{equation}
It immediately follows that, for a given configuration, the computational complexity: 1) increases proportionally to the square of $r^\textrm{p-c}$ when fixing the parent mesh number $M^\textrm{p}$, and 2) decreases with the square of $r^\textrm{p-c}$ when fixing the child mesh number $M^\textrm{c}$. These remarkable characteristics are valid only for $\text{M}^2$-FEM and are not applicable to traditional nonlocal finite element approaches because, the relative resolution is fixed as $r^\textrm{p-c}=1$ following the direct sharing of the mesh between the parent and child domains. To complement the above theoretical deduction, we present performance data from simulations in the Fig.~(\ref{fig: computational_analysis}).

Figure~(\ref{fig: computational_analysis}) compares the variation of the computational cost (in terms of the average time needed for a set of simulations) against the relative error (in terms of the difference in simulated $\textrm{max}(u_x)$ from the reference value) for Case Study 1. All the simulations were performed on a personal computer equipped with an Intel(R) Core(TM) i7-9750H processor and 16 GB RAM. Figure~(\ref{fig: computational_analysis}-a) presents the growth of total computational time with respect to both $M^\textrm{p}$ and $r^\textrm{p-c}$. As evident from the expression of $M^\textrm{sum}$ in Eq.~(\ref{eq: M_sum}), the time complexity can be controlled by varying either the parent ($M^\textrm{p}$) or child mesh size ($r^\textrm{p-c}$). Figure~(\ref{fig: computational_analysis}-b) presents the variation of time cost and error with respect to $M^\textrm{p}$ and $r^\textrm{p-c}$. As evident, by tuning the values of $M^\textrm{p}$ and $r^\textrm{p-c}$, $\text{M}^2$-FEM can be highly efficient (reflect in its ability to save computational time) while still maintaining the same level of accuracy. This aspect is more transparent from two sets of simulations in the Fig.~(\ref{fig: computational_analysis}):
\begin{itemize}[leftmargin=*]
    \item the solution obtained at $(M^\textrm{p},r^\textrm{p-c})=(10,2)$ requires (roughly) the same simulation time as that at $(M^\textrm{p},r^\textrm{p-c})=(15,1)$, but has much smaller error (see the data highlighted by dotted cyan circles),
    
     \item the solution obtained at $(M^\textrm{p},r^\textrm{p-c})=(10,3)$ has roughly the same error as that at $(M^\textrm{p},r^\textrm{p-c})=(25,1)$, but requires less computational time (see the data highlighted by dotted magenta circles).
\end{itemize}
Notably, for a fixed parent mesh, the simulation using a finer child mesh ($r^\textrm{p-c}>1$) shows better performance than the one with a non-refined child mesh ($r^\textrm{p-c}=1$). These results indicate that, when compared to existing (mesh-coupling-based) integral solvers, $\text{M}^2$-FEM equipped with mesh-decoupling technique admits greater computational flexibility; ultimately, leading to higher accuracy and higher efficiency within fixed computational resources.

\section{Summary and outlook}\label{sec: 6}
This study presented the formulation, the numerical implementation, and the performance of $\text{M}^2$-FEM that is an advanced finite element solver capable of simulating generalized integral theories. $\text{M}^2$-FEM derives its unique characteristics from the concept of multiscale mesh-decoupling and scale-bridging, which combined create a numerically independent and localized set of meshes to approximate the underlying integral theory at different material scales, in a physically consistent manner. The characteristic multiscale independence has a multifold impact on the overall computational solver. More specifically, it enables the independent topology-tuning of the scale-specific localized meshes such that $\text{M}^2$-FEM is able to simulate complex multiscale configurations and admit generalized kernel functions. These latter characteristics are achieved while simultaneously providing high accuracy and computational efficiency. These features enable $\text{M}^2$-FEM to model the static and dynamic behavior of a wide class of advanced materials that exhibit unconventional geometry, nonlocal and multiscale effects, and anomalous dispersion behavior. 

Overall, $\text{M}^2$-FEM removes significant computational drawbacks in the application of integral structural theories. While the present study was developed based on Eringen's strain-driven nonlocal theory, $\text{M}^2$-FEM is very general in nature and can be readily applied to other integral theories including, for example, self-consistent approaches, homogenization approaches, multiscale averaging approaches, and displacement-driven nonlocal approaches. The extension of $\text{M}^2$-FEM to these approaches merely requires a different numerical evaluation of the child-domain integral contributions by following the model-specific constitutive theory. The underlying mesh-decoupling and scale-bridging, fundamental pillars of $\text{M}^2$-FEM, remain unaltered. The recursive nesting ability of $\text{M}^2$-FEM highlighted the ability of the proposed method to support the analysis of advanced material models with hierarchical nesting. These characteristics suggest that $\text{M}^2$-FEM can serve as the foundation for a next-generation class of solvers offering critical capabilities to enable the simulation of highly multiscale systems.

\paragraph*{Data Availability.}~All the necessary data and information required to reproduce the results are available in the paper and the supplementary information document.

\paragraph*{Declaration of Competing Interest.}~The authors declare that they have no known competing financial interests or personal relationships that could have appeared to influence the work reported in this paper.

\paragraph*{Acknowledgements.}~The authors gratefully acknowledge the financial support of the National Science Foundation under grants MOMS \#1761423 and DCSD \#1825837. S.P. acknowledges the partial support received through Dimitris N. Chorafas Foundation Award (2022). Any opinions, findings, and conclusions or recommendations expressed in this material are those of the author(s) and do not necessarily reflect the views of the National Science Foundation.

\bibliographystyle{naturemag}
{{\bibliography{Report}}}
\end{document}